\documentclass[10pt]{amsart}
\usepackage{amssymb}
\usepackage{dsfont}
\usepackage{amscd}
\usepackage[mathscr]{euscript} 
\usepackage[all]{xy}
\usepackage{url}
\usepackage{stmaryrd}
\usepackage{comment}
\usepackage[retainorgcmds]{IEEEtrantools}
\usepackage{hyperref}

\title{PAC structures as invariants of finite group actions}
\author[D. M. HOFFMANN]{Daniel Max Hoffmann$^{\dagger}$}
\thanks{$^{\dagger}$SDG. Supported by the Narodowe Centrum Nauki grant no. 2021/43/B/ST1/00405.}
\address{$^{\dagger}$ Instytut Matematyki\\
Uniwersytet Warszawski\\
Warszawa\\
Poland}
\email{daniel.max.hoffmann@gmail.com}
\urladdr{https://sites.google.com/site/danielmaxhoffmann/home}
\author[P. KOWALSKI]{Piotr Kowalski$^{\ddagger}$}
\thanks{$^{\ddagger}$ 
 Supported by the Narodowe Centrum Nauki grant no. 2021/43/B/ST1/00405 and by the T\"{u}bitak 1001 grant no. 119F397.}
\address{$^{\ddagger}$Instytut Matematyczny\\
Uniwersytet Wroc{\l}awski\\
Wroc{\l}aw\\
Poland}
\email{pkowa@math.uni.wroc.pl} \urladdr{http://www.math.uni.wroc.pl/\textasciitilde pkowa/ }

\DeclareMathOperator{\locus}{locus}\DeclareMathOperator{\cl}{cl}
\DeclareMathOperator{\acl}{acl} \DeclareMathOperator{\dcl}{dcl} 
 \DeclareMathOperator{\aut}{Aut} \DeclareMathOperator{\id}{id}

 \DeclareMathOperator{\dom}{dom}

 \DeclareMathOperator{\theo}{Th}
 \DeclareMathOperator{\eq}{eq}
 
\DeclareMathOperator{\tp}{tp}

\DeclareMathOperator{\rat}{rat}
\DeclareMathOperator{\sep}{sep}

\DeclareMathOperator{\dcf}{DCF}\DeclareMathOperator{\scf}{SCF}
\DeclareMathOperator{\acf}{ACF}

\DeclareMathSymbol{\shortminus}{\mathbin}{AMSa}{"39}

\DeclareMathOperator{\FC}{\mathfrak{C}}

\newtheorem{theorem}{Theorem}[section]
\newtheorem{prop}[theorem]{Proposition}
\newtheorem{lemma}[theorem]{Lemma}
\newtheorem{cor}[theorem]{Corollary}
\newtheorem{fact}[theorem]{Fact}

\theoremstyle{definition}
\newtheorem{definition}[theorem]{Definition}
\newtheorem{example}[theorem]{Example}
\newtheorem{remark}[theorem]{Remark}
\newtheorem{question}[theorem]{Question}

\newtheorem{innercustomthm}{Theorem}
\newenvironment{customthm}[1]
  {\renewcommand\theinnercustomthm{#1}\innercustomthm}
  {\endinnercustomthm}

\begin{document}
\newcommand{\fg}{\mathfrak{g}}

\newcommand{\lili}{\underleftarrow{\lim }}
\newcommand{\coco}{\underrightarrow{\lim }}
\newcommand{\twoc}[3]{ {#1} \choose {{#2}|{#3}}}
\newcommand{\thrc}[4]{ {#1} \choose {{#2}|{#3}|{#4}}}
\newcommand{\Zz}{{\mathds{Z}}}
\newcommand{\Ff}{{\mathds{F}}}
\newcommand{\Cc}{{\mathds{C}}}
\newcommand{\Rr}{{\mathds{R}}}
\newcommand{\Nn}{{\mathds{N}}}
\newcommand{\Qq}{{\mathds{Q}}}
\newcommand{\Kk}{{\mathds{K}}}
\newcommand{\Pp}{{\mathds{P}}}
\newcommand{\ddd}{\mathrm{d}}
\newcommand{\Aa}{\mathds{A}}
\newcommand{\dlog}{\mathrm{ld}}
\newcommand{\ga}{\mathbb{G}_{\rm{a}}}
\newcommand{\gm}{\mathbb{G}_{\rm{m}}}
\newcommand{\gaf}{\widehat{\mathbb{G}}_{\rm{a}}}
\newcommand{\gmf}{\widehat{\mathbb{G}}_{\rm{m}}}
\newcommand{\ka}{{\bf k}}
\newcommand{\ot}{\otimes}
\newcommand{\si}{\mbox{$\sigma$}}
\newcommand{\ks}{\mbox{$({\bf k},\sigma)$}}
\newcommand{\kg}{\mbox{${\bf k}[G]$}}
\newcommand{\ksg}{\mbox{$({\bf k}[G],\sigma)$}}
\newcommand{\ksgs}{\mbox{${\bf k}[G,\sigma_G]$}}
\newcommand{\cks}{\mbox{$\mathrm{Mod}_{({A},\sigma_A)}$}}
\newcommand{\ckg}{\mbox{$\mathrm{Mod}_{{\bf k}[G]}$}}
\newcommand{\cksg}{\mbox{$\mathrm{Mod}_{({A}[G],\sigma_A)}$}}
\newcommand{\cksgs}{\mbox{$\mathrm{Mod}_{({A}[G],\sigma_G)}$}}
\newcommand{\crats}{\mbox{$\mathrm{Mod}^{\rat}_{(\mathbf{G},\sigma_{\mathbf{G}})}$}}
\newcommand{\crat}{\mbox{$\mathrm{Mod}^{\rat}_{\mathbf{G}}$}}
\newcommand{\cratinv}{\mbox{$\mathrm{Mod}^{\rat}_{\mathbb{G}}$}}
\newcommand{\ra}{\longrightarrow}
\newcommand{\bdcf}{B-\dcf}
\makeatletter
\providecommand*{\cupdot}{%
  \mathbin{%
    \mathpalette\@cupdot{}%
  }%
}
\newcommand*{\@cupdot}[2]{%
  \ooalign{%
    $\m@th#1\cup$\cr
    \sbox0{$#1\cup$}%
    \dimen@=\ht0 %
    \sbox0{$\m@th#1\cdot$}%
    \advance\dimen@ by -\ht0 %
    \dimen@=.5\dimen@
    \hidewidth\raise\dimen@\box0\hidewidth
  }%
}

\providecommand*{\bigcupdot}{%
  \mathop{%
    \vphantom{\bigcup}%
    \mathpalette\@bigcupdot{}%
  }%
}
\newcommand*{\@bigcupdot}[2]{%
  \ooalign{%
    $\m@th#1\bigcup$\cr
    \sbox0{$#1\bigcup$}%
    \dimen@=\ht0 %
    \advance\dimen@ by -\dp0 %
    \sbox0{\scalebox{2}{$\m@th#1\cdot$}}%
    \advance\dimen@ by -\ht0 %
    \dimen@=.5\dimen@
    \hidewidth\raise\dimen@\box0\hidewidth
  }%
}
\makeatother

\def\Ind#1#2{#1\setbox0=\hbox{$#1x$}\kern\wd0\hbox to 0pt{\hss$#1\mid$\hss}
\lower.9\ht0\hbox to 0pt{\hss$#1\smile$\hss}\kern\wd0}

\def\ind{\mathop{\mathpalette\Ind{}}}

\def\notind#1#2{#1\setbox0=\hbox{$#1x$}\kern\wd0
\hbox to 0pt{\mathchardef\nn=12854\hss$#1\nn$\kern1.4\wd0\hss}
\hbox to 0pt{\hss$#1\mid$\hss}\lower.9\ht0 \hbox to 0pt{\hss$#1\smile$\hss}\kern\wd0}

\def\nind{\mathop{\mathpalette\notind{}}}

\maketitle

\begin{abstract}
We study model theory of actions of finite groups on substructures of a stable structure.
We give an abstract description of existentially closed actions as above in terms
of invariants and PAC structures.
We show that if the corresponding PAC property is first order, then the theory of such actions has a model companion. Then, we analyze some particular theories of interest (mostly various theories of fields of positive characteristic) and show that in all the cases considered the PAC property is first order.
\end{abstract}

\section{Introduction}

In this paper, we consider the notion of a \emph{pseudo algebraically closed (PAC) substructure} of a stable structure. This notion originates from the theory of pseudo algebraically closed fields, which were first considered by Ax in 1960's while he worked on pseudofinite fields (\cite{ax68}). Studying PAC structures beyond the case of fields was initiated by Hrushovski (\cite{manuscript}) in the strongly minimal context. Pillay and Polkowska considered the PAC property in the stable case (\cite{PilPol}), there are slight differences with the approach we take here.
 PAC structures also appeared in Afshordel's thesis (\cite{Afshordel}). Recently, PAC structures were analized by the first author (\cite{Hoff3}, \cite{Hoff4}) and also by Dobrowolski, the first author, and Lee (\cite{DHL1}).


Here, we are working with a (complete) stable theory $T$ which admits quantifier elimination and then focus on its universal part $T_{\forall}$. In other words, a typical situation looks as follows.
We have a universal theory $T_{\forall}$ with a stable model completion $T$, so $T$ has quantifier elimination and $T$ axiomatizes existentially closed models of $T_{\forall}$.
Then, intuitively, the class of PAC structures in $T$ lies in between the class of existentially closed structures (models of $T$) and the class of all the structures considered (models of $T_{\forall}$). There are several possible definitions of the notion of PAC, we adopt here the definition from \cite{Hoff3} (expressed in terms involving stationary types), which is a slight modification of the definition from \cite{PilPol}, and which is equivalent to Afshordel's definition from \cite{Afshordel} in the case of stable theories.
To define the notion of a PAC structure, one needs to use an appropriate notion of irreducibility. In the classical case of PAC fields, a topological notion is used coming from the Zariski topology. Hrushovski used in \cite{manuscript} ``Morley irreducibility'', that is he considered definable sets of Morely degree one. Pillay and Polkowska used \cite{PilPol} stationary types and we proceed similarly here (however, we avoid any saturation requirements as given in \cite{PilPol}).
We say that a structure $F\models T_{\forall}$ is PAC in $T$ (see Definition \ref{def:PAC}) if all stationary types (in the sense of the theory $T$) over $F$ are finitely satisfiable in $F$. Let us point out that in the case of the theory of algebraically closed fields, all the irreducibility notions mentioned above are essentially the same. However, this is not the case for other theories of interest as the theory of differentially closed fields of characteristic $0$ or the theory of compact complex manifolds (see Section \ref{secnoeth}). Nevertheless, we show in Section \ref{secgen} that all these irreducibility notions lead to the same notion of a PAC structure.

For an extension $F\subseteq K$ of models of $T_{\forall}$, we obtain relative notions of \emph{$K$-strongly PAC} and \emph{algebraically $K$-strongly PAC} (see Definition \ref{def:PAC2}). They are meaningful and can be though of as measuring the distance between being PAC and being a model of $T$ ($K$-strongly PAC) or between being definably closed and algebraically closed (algebraically $K$-strongly PAC), see Remark \ref{newrem}.

Our main motivation for considering PAC structures comes from model theory of group actions. In the set-up above, we consider actions of a fixed group $G$ on models of $T_{\forall}$ by automorphisms. Clearly, such actions are first-order expressible in an appropriate language and we aim to describe existentially closed actions and check whether a model companion of the theory of such actions exists. The result below may be considered as an abstract generalization of our theorem about finite group actions on fields (see \cite[Theorem 3.29]{HK3})
and as a continuation of studies from \cite{Hoff3}.

\begin{customthm}{\ref{thm.4conditions}}
Let $G$ be a finite group and $T$ be a stable theory coding finite sets, which has quantifier elimination and eliminates strong types (that is: types over algebraically closed sets are stationary). Assume that $G$ acts faithfully on $K=\dcl(K)\models T_{\forall}$. Then, the following are equivalent.
\begin{enumerate}
  \item The action of $G$ on $K$ is existentially closed.

  \item The structure of invariants $K^G$ is $K$-strongly PAC.

  \item The structure of invariants $K^G$ is PAC and algebraically $K$-strongly PAC.
\end{enumerate}
\end{customthm}
The above theorem gives a description of existentially closed finite group actions, but it is not clear whether this description is first-order, so this theorem does not settle the question of the existence of a model companion of the theory of finite actions. We can show the following implication.

\begin{customthm}{\ref{thm:ec.if.PAC}}
Let $G$ be a finite group and $T$ be as in the statement of Theorem \ref{thm.4conditions}. If the class of $T$-PAC structures is elementary, then the model companion of the theory of $G$-actions on models of $T_{\forall}$ exists.
\end{customthm}
After the abstract description of existentially closed actions (Theorem \ref{thm.4conditions}) and giving a criterion for existence of a model companion of the theory of finite actions (Theorem \ref{thm:ec.if.PAC}), we focus on particular examples of theories. We discuss the following three stable theories of fields of positive characteristic ($p$ is a prime and $e$ is a positive integer):
\begin{enumerate}
\item The theory $\mathrm{SCF}_{p,e}$ of separably closed fields of characteristic $p$ and inseparability degree $e$.

\item The theory $\mathrm{SCF}_{p,\infty}$ of separably closed fields of characteristic $p$ and infinite inseparability degree.

\item The theory $\mathrm{DCF}_{p}$ of differentially closed fields of characteristic $p$.
\end{enumerate}
In the most interesting cases of the theories $\mathrm{SCF}_{p,\infty}$ and $\mathrm{DCF}_{p}$, we do not have elimination of imaginaries, however we still have its weaker versions (coding finite sets and eliminating strong types), which are enough for the set-up from Theorems \ref{thm.4conditions} and \ref{thm:ec.if.PAC}. For these theories, we describe PAC structures in a first-order way using a result of Tamagawa (see Theorem \ref{tamagawa}) about positive characteristic PAC fields. We finish with some general questions regarding the PAC property and existence of a model companion of the theory of finite actions. It should be mentioned that after replacing a finite group $G$ with the infinite cyclic group $(\Zz,+)$, then the model theory of actions of $(\Zz,+)$ has been thoroughly studied (see e.g. \cite{acfa1} and \cite{ChPi}). We compare these two situation in Section \ref{secoeq}.

This paper is organized as follows. In Section 2, we introduce several versions of the notion of a PAC structure and show the basic results about them. In Section 3, we put the group action to the picture and prove the main two abstract results stated above (Theorems \ref{thm.4conditions} and \ref{thm:ec.if.PAC}). In Section 4, we consider some particular theories (mostly theories of fields of positive characteristic) and give a first order characterization of PAC structures with respect to these theories.

\section{Preliminaries}
\subsection{Set-up}
Let $T$ be a complete first order theory with
a monster model $\FC\models T$ (i.e. a strongly $\bar{\kappa}$-homogeneuos and $\bar{\kappa}$-saturated model of $T$ for a very big cardinal $\bar{\kappa}$).
Throughout the paper, $\acl$ and $\dcl$ mean the
algebraic closure and the definable closure in $\FC$. Usually, $x$ stands for a (finite) tuple of variables.
Moreover, for the rest of this paper, let $G$ be a group such that $|G|<\bar{\kappa}$.

Bearing in mind any future applications, 
we try in this paper to formulate each result with a minimal list of assumptions.
Therefore, we organize our general model-theoretic assumptions in the following list (we are aware that there are some overlaps, but we preferred more transparent exposition):
\begin{enumerate}

\item[(QE)] $T$ has quantifier elimination.

\item[(FS)] $T$ codes finite tuples (i.e. eliminates finite imaginaries).

\item[(FS$+$)] $T$ has (FS) and for every $k<\omega$, for every variable $x$ corresponding to a real sort and the $0$-definable equivalence relation $E$ on $S_x^k$ given by 
$$E(\bar{x},\bar{x}')\quad\iff\quad\{x_1,\ldots,x_k\}=\{x_1',\ldots,x_k'\},$$
there exists a $0$-definable in $L$ function $f:S_x^k\to S_w$ such that $E$ is a fibration of $f$.


\item[(ST$+$)] $T$ is stable and types over algebraically closed sets are stationary (elimination of strong types).
\end{enumerate}
Convention: if a statement starts with any combination of the above properties, it means that we assume the properties given in this particular combination. For example, the following remark assumes property (FS):

\begin{remark}
(FS)
The condition (FS$+$) is equivalent to:
\begin{itemize}
\item
on each sort there is at least one $0$-definable element, and
\item 
there is a sort with at least two $0$-definable elements.
\end{itemize}
\end{remark}

\begin{proof}
Similarly as in the proof of Lemma 8.4.7 from \cite{tentzieg}, but, here, we allow many sorted structures.
\end{proof}

\begin{remark}
Let us discuss what one can do to meet the above requirements if starting from arbitrary stable $L_0$-theory $T_0$.
As we would like to work under assumptions of quantifier elimination and elimination of imaginaries,
we pass to the language $L:=(L_0^{\eq})^m$ and $L$-theory $T:=(T_0^{\eq})^m$ (we add imaginary sorts and then do the Morleysation). This new theory $T$ is stable, has quantifier elimination and elimination of imaginaries.
On top of that, 
every $0$-definable equivalence relation $E$ on $\FC^n$ is the fibration of 
the canonical projection $\pi_E:\FC^n\to\FC^n/E$ which is build-in in the language $(L_0^{\eq})^m$, thus
a $0$-definable function. Strong types in any stable theory are stationary, and $b\ind_A A$ for any $b$ and $A$. Therefore $T$ enjoys all the properties: (QE), (FS), (FS$+$)
and (ST$+$).
\end{remark}

\subsection{Notion of PAC structure and auxiliary facts}
In this subsection, we recall
several definitions and useful facts from \cite{Hoff3} and \cite{Hoff4}. We also provide a few new notions closely related to the old definitions.
The reader may also consult \cite{PilPol} and \cite{Polkowska} for more on PAC structures in general model theoretic framework. Also \cite{Afshordel} provides a nice of exposition of the notion of a PAC structure and related topics. A well-written survey on different variants of the notion of elimination of imaginaries and related concepts from the Galois theory is \cite{casfar}.

\begin{definition}\label{def:PAC}
(Let $T$ be stable.)
A substructure $F$ of $\FC$ is \emph{pseudo-algebraically closed} (\emph{PAC}) if every stationary type over $F$ (in the sense of the $L(F)$-theory of $\FC$) is finitely satisfiable in $F$.
\end{definition}

\noindent
The above definition appears in \cite{Hoff3} (see also Definition 5.29 in \cite{Afshordel}). In subsection 3.1 of \cite{Hoff3}, there is a discussion on 
possible choices of the definition
of a PAC substructure and a comparison of Definition \ref{def:PAC} to definitions of PAC structures given in \cite{manuscript} and in \cite{PilPol}. In short, Definition \ref{def:PAC} coincides with the definition of a PAC substructure in the strongly minimal context of \cite{manuscript} and relaxes the saturation assumption from the definition of a PAC substructure from \cite{PilPol}.
Note that every PAC substructure is automatically definably closed. Thus PAC substructures for $T=\acf$ coincide with \textbf{perfect} pseudo-algebraically closed fields (as defined in e.g. \cite{FrJa}).

\begin{definition}\label{def:PAC2}
Let $F=\dcl(F)\subseteq K\subseteq \FC$. 
\begin{enumerate}
    \item 
    We say that $F$ is \emph{$K$-strongly PAC} if each type $p(x)\in S(F)$, which has a unique non-forking extension over $K$, is finitely satisfiable in $F$.
    
    \item
    We say that $F$ is \emph{algebraically $K$-strongly PAC} if each algebraic type $p(x)\in S(F)$, which has a unique non-forking extension over $K$, is finitely satisfiable (thus realized) in $F$.
\end{enumerate}
\end{definition}

Note that being $K$-strongly PAC for $F\subseteq K$ implies being 
algebraically $K$-strongly PAC for $F$.
Moreover, being $K$-strongly PAC for $F$ implies being a PAC substructure for $F$.

\begin{remark}\label{newrem}
It should help to understand the relative notions of (algebraically) $K$-strongly PAC by considering the ultimate cases of $K=F$ and $K\models T$. It is quite easy to see the following.
\begin{enumerate}
  \item A structure $F$ is $F$-strongly PAC if and only if $F\models T$.
  \item ($T$ is stable) A structure $F$ is $K$-strongly PAC for $K\models T$ if and only if $F$ is PAC.
  \item A structure $F$ is algebraically $F$-strongly PAC if and only if $F=\acl(F)$.
  \item A structure $F$ is algebraically $K$-strongly PAC for $K\models T$ if and only if $F=\dcl(F)$.
  
\end{enumerate}
\end{remark}

\begin{definition}\label{regular.def}
\begin{enumerate}
\item Let $F\subseteq K$ be small subsets of $\mathfrak{C}$. We say that $F\subseteq K$ is \emph{primary} if
$$\dcl(K)\cap\acl(F)=\dcl(F).$$

\item
Let $F\subseteq K$ be small subsets of $\mathfrak{C}$.
We say that $F\subseteq K$ is \emph{regular} if $F\subseteq K$ is primary and $F=\dcl(F)$.

\item Let $F$ be a small definably closed substructure of $\mathfrak{C}$. 
We say that $F$ is \emph{regularly closed}
if for every small substructure $F'$ of $\mathfrak{C}$, which is a regular extension of $F$,
it follows $F\preceq_1 F'$ (i.e. $F$ is existentially closed in $F'$).
\end{enumerate}
\end{definition}

\noindent
The above notion of a primary extension was previously (e.g. \cite{Hoff3}, \cite{Hoff4}) called ``regular''. It corresponds to regular extensions in $T=$ACF provided the smaller field is perfect (equivalently, definably closed). Here, we decided to follow closer the terminology from the theory of fields and distinguish between ``primary'' and ``regular'' extensions.
We plan to refine even more the notion of the model-theoretic ``regular'' extension after studying a possible notion of the model-theoretic separable extension in the future.

Now, we will sharpen facts from earlier articles which lead to the main results in this manuscript.
The majority of \cite{Hoff3} was written under the assumption of (full) elimination of imaginaries, elimination of quantifiers and stability.
This is fine if we are interested in an abstract approach to the subject. 
However, as we are interested in applications of our results to particular theories, which do not enjoy elimination of imaginaries (see Section \ref{sec:examples}), we need to relax this assumption. Moreover, the assumption on stability was not crucial in several useful facts from \cite{Hoff3}, making them applicable in a broader context.
Therefore we take the opportunity to provide the following results with minimal assumptions.
The proofs of the following facts remain almost the same as
the proofs of
their counterparts from \cite{Hoff3}. Recall that ``regular'' extensions from \cite{Hoff3}
are now ``primary'' extensions.

All the proper subsets, substructures and tuples of the monster model $\FC$ are, if not stated otherwise, small in comparison to the saturation of $\FC$. Here, upper case letters, like $E$ or $A$, are denoting proper subsets, and lower case letters, like $a$, stand for tuples.

\begin{fact}[Fact 3.32 in \cite{Hoff3}]\label{fact.3.32}
(FS) If $E\subseteq A$ is primary then for every $a\in\acl(E)$ there exists a unique extension of $\tp(a/E)$ over $A$.
\end{fact}

\begin{fact}[Fact 3.33 in \cite{Hoff3}]\label{fact.3.33}
(FS) If $E\subseteq A$ is primary , $f_1,f_2\in\aut(\FC)$ and $f_1|_E=f_2|_E$,
then there exists $h\in\aut(\FC)$ such that $h|_A=f_1|_A$ and $h|_{\acl(E)}=f_2|_{\acl(E)}$.
\end{fact}

\begin{fact}[Corollary 3.34 in \cite{Hoff3}]\label{fact.3.34}
(FS) If $E\subseteq A$ is primary and $A_0\subseteq A$ then $\tp(A_0/E)$ has a unique extension over $\acl(E)$.
\end{fact}

The following definition is taken from page 21. of \cite{Afshordel}.

\begin{definition}\label{def:acl.stationary}
We say that a type $p(x)\in S(A)$ is $\acl$-stationary if it has a unique extension over $\acl(A)$.
\end{definition}

\begin{lemma}
(FS) Consider $p\in S(E)$. The following are equivalent:
\begin{enumerate}
\item $p$ is $\acl$-stationary,
\item $E\subseteq \dcl(Ea)$ is primary for some $a\models p$,
\item $E\subseteq \dcl(Ea)$ is primary for every $a\models p$.
\end{enumerate}
\end{lemma}

\begin{proof}
The proof is similar to the proof of Lemma 3.35 in \cite{Hoff3}, but a few steps require sharper reasoning, thus we include it here.

The equivalence (2)$\iff$(3) follows by definition.
First, we argue for (1)$\Rightarrow$(2): assume (1) and suppose that (2) does not hold.
As $p$ is $\acl$-stationary, there exists a unique extension 
$p|_{\acl(E)}$ of $p$ over $E$. Let $a\models p|_{\acl(E)}$, then $a\models p$ and $E\subseteq Ea$ is not primary.
Take 
$$c\in\dcl(Ea)\cap\acl(E)\setminus\dcl(E).$$
Since $c\not\in\dcl(E)$, there exists $f\in\aut(\FC/E)$ such that $f(c)\neq c$.
We see that $f(a)\models p|_{\acl(E)}$, so there exists $h\in\aut(\FC/\acl(E))$ such that $h(a)=f(a)$.
Note that $h^{-1}f\in\aut(\FC/Ea)$ and, because $c\in\dcl(Ea)$ and $c\in\acl(E)$,
$$c=h^{-1}f(c)=f(c)\neq c,$$
so a contradiction.
The implication (2)$\Rightarrow$(1) is contained in Fact \ref{fact.3.34}.
\end{proof}

\begin{fact}[Lemma 3.35 in \cite{Hoff3}]\label{fact.3.35}
(FS, ST$+$) Consider $p\in S(E)$.
The following are equivalent:
\begin{enumerate}
\item $p$ is stationary,
\item $p$ is $\acl$-stationary,
\item $E\subseteq Ea$ is primary for some $a\models p$,
\item $E\subseteq Ea$ is primary for every $a\models p$.
\end{enumerate}
\end{fact}

\begin{fact}[Corollary 3.36 in \cite{Hoff3}]\label{fact.3.36}
(QE, FS, ST$+$) For any small substructure $N$ there exists a non-algebraic stationary type over $N$ in any finitely many variables.
\end{fact}

\begin{fact}[Corollary 3.38 in \cite{Hoff3}]\label{fact.3.38}
(FS, ST$+$) Assume that $A,B\subseteq\FC$, $E\subseteq A$ is primary, $f_1,f_2\in\aut(\FC)$ and $f_1|_E=f_2|_E$.
If $A\ind_E B$ and $f_1(A)\ind_{f_1(E)} f_2(B)$ then there exists $h\in\aut(\FC)$ such that $h|_A=f_1|_A$ and $h|_B=f_2|_B$.
\end{fact}

\begin{fact}[Lemma 3.39 in \cite{Hoff3}]\label{fact.3.39}
(FS, ST$+$) If $E\subseteq A\cap B$, $E\subseteq A$ is primary and $B\ind_E A$ then $B\subseteq BA$ is primary.
\end{fact}

\begin{fact}[Corollary 3.40 in \cite{Hoff3}]\label{fact.3.40}
(FS, ST$+$) If $E\subseteq A$ and $E\subseteq B$ are primary, and $B\ind_E A$ then also $E\subseteq BA$ is primary.
\end{fact}

\begin{remark}
\begin{enumerate}
\item
(FS, ST$+$) $F\subseteq K$ is primary if and only if for every tuple $b$ from $\dcl(K)$, the type $\tp(b/F)$ is stationary (Fact \ref{fact.3.35}).

\item
(QE, FS, ST$+$) Using the item (1), a substructure $F$ is PAC if and only if it is definably closed and regularly closed.
\end{enumerate}
\end{remark}

\begin{definition}\label{galois.ext.def}
\begin{enumerate}
\item Assume that $F\subseteq K$ are substructures of $\FC$. 
We say that $K$ is \emph{normal over $F$} (or we say that $F\subseteq K$ is a \emph{normal extension}) if 
$\sigma(K)\subseteq K$ for every $\sigma\in\aut(\FC/K)$.
(Note that if $K$ is small and $F\subseteq K$ is normal, then it must be $K\subseteq\acl(F)$.)

\item Assume that $F\subseteq K\subseteq\acl(F)$ are small substructures of $\mathfrak{C}$ 
such that $F=\dcl(F)$, $K=\dcl(K)$ and $K$ is normal over $F$. In this situation we say that $F\subseteq K$ is a \emph{Galois extension}.
\end{enumerate}
\end{definition}

\begin{definition}
Assume that $F\subseteq K$ is an extension of substructures in $\FC$.
We define the Galois group of the extension $F\subseteq K$ as
$$G(K/F):=\aut(K/F)=\{f|_K\;|\;f\in\aut(\FC/F),\,f(K)=K \}.$$
Moreover $B$ is any subset of $\FC$, then the extension $\dcl(B)\subseteq\acl(B)$ is Galois and we speak about the absolute Galois group of $B$ which is the following profinite group:
$$G(B):=G(\acl(B)/\dcl(B)).$$
\end{definition}

Note that the above definition of $G(K/F)$ is often expressed in terms of the automorphisms of $K$ as an $L$-structure on its own, but as we will work under the assumption of the quantifier elimination, both variants of the definition coincide and it just the matter of taste.

The following useful fact is standard and its proof is straightforward.

\begin{lemma}\label{lemma.transitively}
Assume that $F\subseteq K$ is a Galois extension and $p(x)\in S(F)$.
Then the Galois group $G(K/F)$ acts transitively on the set of extensions of $p$ over $K$.
\end{lemma}

The following definition and example are taken from \cite{DHL1} and \cite{PilPol}.
A more detailed discussion of examples of PAC structures and the property from Definition \ref{def:PAC.first.order} will be given in Section \ref{sec:examples}.

\begin{definition}\label{def:PAC.first.order}
(Let $T$ be stable.)
We say that \emph{PAC is a first order property in $T$} ($=\theo(\FC)$) if there exists a set $\Sigma$ of $\mathcal{L}$-sentences such that for any $P\subseteq \FC$
$$P\models \Sigma\qquad\iff\qquad P\text{ is PAC}.$$
\end{definition}

\begin{example}
\begin{enumerate}
\item
PAC is a first order property in ACF$_p$ for $p=0$ and for $p$ being a prime number, see Proposition 11.3.2 in \cite{FrJa}.

\item 
The axioms given in Proposition 5.6 from \cite{PilPol} show that 
PAC is a first order property (in the above sense) in DCF$_0$
which is formulated in a different way than
the condition ``PAC is a first order property" appearing in \cite{PilPol}. 
\end{enumerate}
\end{example}

\section{Finite group actions}
The main goal of this section is to describe existentially closed substructures with a finite group action in first order terms. The general strategy is as follows. First, characterize their structure by the structure of the invariants of the group action, then answer which properties of the invariants correspond to the existential closedeness of the whole substructure with group action. Finally, express these properties as first order statements.

\subsection{Basic facts}
We introduce the language $L_G$ being the language $L$ extended by a unary function symbol $\sigma_g$ for each $g\in G$, i.e. $L_G=L\,\cup\,\{\sigma_g\;|\;g\in G\}$. Often, ``$\sigma_g$'' will denote also the interpretation of the symbol $\sigma_g$ in a given $L_G$-structure.
Moreover, we set $\bar{\sigma}:=(\sigma_g)_{g\in G}$.
We consider the collection of sentences in the language $L_G$, say $A_G$, which precisely expresses the following
\begin{itemize}
\item $\sigma_g$ is an automorphism of the $L$-structure for every $g\in G$,
\item $\sigma_g\circ\sigma_h=\sigma_{g\cdot h}$ for all $g,h\in G$.
\end{itemize}
In other words, if $K$ is an $L$-structure, and there exists an $L_G$-structure $(K,\bar{\sigma})$ living on $K$, 
we have that $(K,\bar{\sigma})\models A_G$ if and only if
for each $g\in G$ we have that $\sigma_g\in\aut(K)$ and the map
$$G\ni g\mapsto\sigma_g\in\aut(K)$$
is a group homomorphism.

\begin{definition}
\begin{enumerate}
\item Let $(K,\bar{\sigma})$ be an $L_G$-structure. We say that $\bar{\sigma}$ is a $G$-action on $K$ if $(K,\bar{\sigma})\models A_G$.

\item If $T'$ is an $L$-theory, then by $(T')_G$ we denote the set of consequences of $T'\,\cup\,A_G$.

\item
If $(K,\bar{\sigma})\models (T_{\forall})_G$, where $K$ is of cardinality smaller than the saturation of $\FC$, then we call it a \emph{substructure with $G$-action}. 
Note that, without loss of generality, $K\subseteq\FC$, thus the name ``substructure''.

\item 
We say that  a substructure with $G$-action
$(K,\bar{\sigma})$ is \emph{existentially closed}
if 
$(K,\bar{\sigma})$ is
an existentially closed model of the theory $(T_{\forall})_G$.

\item 
If the existentially closed models of the theory $(T_{\forall})_G$ form an elementary class, we denote the theory of this class by $G-T$.

\end{enumerate}
\end{definition}

\begin{definition}
Assume that $(K,\bar{\sigma})$ is a substructure with $G$-action.
Then we denote 
$$K^G:=\{a\in K\;|\;(\forall g\in G)\,(\sigma_g(a)=a)\,\}$$
and call it the \emph{substructure of invariants}.
\end{definition}

\begin{remark}\label{rem:dcl}
(QE) Let $(K,\bar{\sigma})$ be a substructure with $G$-action.
If $(K,\bar{\sigma})$ is existentially closed then $K=\dcl(K)$.
If $K=\dcl(K)$ then $K^G=\dcl(K^G)$. For the standard proofs, the reader may consult Remark 3.24 and Remark 3.26 in \cite{Hoff3}.
\end{remark}

\begin{lemma}\label{lemma:extending}
(QE) Let $(K,\bar{\sigma})$ be a substructure with $G$-action and let $p(x)\in S(K)$ be a $G$-invariant type (i.e. $\sigma_g(p)=p$ for every $g\in G$). Then for any $a\models p$ the set $\dcl(K,a)$ might be equipped with a $G$-action extending $(K,\bar{\sigma})$ and acting trivially on $a$.
\end{lemma}

\begin{proof}
Let $a\models p$ and let $\bar{k}$ be some enumeration of $K$.
Then $\bar{k}a\equiv \sigma_g(\bar{k})a$ for any $g\in G$.
This implies that, for each $g\in G$, there exists $\sigma_g'\in\aut(\FC)$ such that $\sigma'_g|_K=\sigma_g$ and $\sigma'_g(a)=a$. Naturally, $(K,(\sigma_g)_{g\in G})\subseteq (\dcl(K,a),(\sigma'_g)_{g\in G})$.
\end{proof}

\begin{fact}[Lemma 2.10 from \cite{Hoff4}]\label{fact:2.10}
(QE, FS) If $G$ is finite and $(K,\bar{\sigma})$ is a substructure with $G$-action such that $\dcl(K)=K$ and the action of $G$ on $K$ is faithful (i.e. if $g\neq h$ then there is $a\in K$ such that $\sigma_g(a)\neq\sigma_h(a)$), then
\begin{itemize}
    \item $K\subseteq \acl(K^G)$,
    \item $K^G\subseteq K$ is a Galois extension,
    \item $G(K/K^G)\cong G$.
\end{itemize}
\end{fact}

\begin{proof}
By Lemma 2.10 from \cite{Hoff4}, Fact 3.7 and Proposition 4.7 from \cite{casfar}.
Being more precise, we obtain the two first bullets as in Lemma 3.23 from \cite{Hoff3}
and then we repeat the proof of Lemma 2.10(4) from \cite{Hoff3} using a variant of the finite Galois correspondence stated in Proposition 4.7 in \cite{casfar}.
\end{proof}


\begin{lemma}\label{lemma:faithful}
(QE, FS, ST$+$) If $(K,\bar{\sigma})$ is an existentially closed substructure with $G$-action, then
the group action if faithful.
\end{lemma}

\begin{proof}
Consider any enumeration of $G$, say $(g_i)_{i\in I}$ where $(I,<)$ is a linear order.
Let $p(x)\in S(K^G)$ be a non-algebraic stationary type (existing by Fact \ref{fact.3.36}), and let $\bar{b}=(b_i)_{i\in I}\models p^{\otimes I}|_{K^G}$ be such that $K\ind_{K^G}\bar{b}$. Let $F$ denote $\dcl(K^G,\bar{b})$, and let $F'$ denote $\dcl(K,\bar{b})$.

As the type $p^{\otimes I}|_K$ is also stationary, the extension $K^G\subseteq F$ is regular.
For each $g\in G$, let $\theta_g$ be a bijection of $I$ such that
$g\cdot g_i=g_{\theta_g(i)}$ holds for each $i\in I$.
As the set $\{b_i\;|\;i\in I\}$ is $K^G$-indiscernible, for each $g\in G$ there exists $\tau_g\in\aut(\FC/K^G)$ such that $\tau_g(b_i)=b_{\theta_g(i)}$.

Now, Corollary 3.38 from \cite{Hoff3}, allows us to simultaneously extend each $\sigma_g$ (over $K$) and $\tau_g$ (over $F$) to an automorphism $\sigma'_g\in\aut(\FC)$, for each $g\in G$.
We have that $(K,(\sigma_g)_{g\in G})\subseteq (F',(\sigma'_g)_{g\in G})$, thus
$(K,(\sigma_g)_{g\in G})$ is existentially closed in  $(F',(\sigma'_g)_{g\in G})$.
If $g\neq h$, then $\sigma'_g(b)\neq\sigma'_h(b)$ for some $b\in F'$,
and so there will be $a\in K$ such that $\sigma_g(a)\neq\sigma_h(a)$.
\end{proof}

\begin{lemma}\label{lemma:Kpac1}
(QE) If $G$ is finitely generated and $(K,\bar{\sigma})$ is an existentially closed substructure with $G$-action, then $K^G$ is $K$-strongly PAC.
\end{lemma}

\begin{proof}
Consider $p(x)\in S(K^G)$ which has a unique non-forking extension over $K$, say $\tilde{p}(x)\in S(K)$.
As $p(x)$ is invariant under action of automorphisms $\sigma_g|_{K^G}$, where $g\in G$, we have that $\tilde{p}(x)$ is invariant under action of automorphisms $\sigma_g$, where $g\in G$ (otherwise, we would get distinct non-forking extensions of $p$ over $K$).

Let $b\models\tilde{p}$, by Lemma \ref{lemma:extending}
there exists an extension of substructures with $G$-action,
$$(K,(\sigma_g)_{g\in G})\subseteq (K',(\sigma'_g)_{g\in G})$$
such that $b\in (K')^G$.
By our assumption, we have that $(K,(\sigma_g)_{g\in G})$ is existentially closed in $(K',(\sigma_g')_{g\in G})$.

Now, let $\varphi(a,x)\in p(x)$. As $T$ has quantifier elimination, we may assume that $\varphi(y,x)$ is quantifier free, what we do.
Of course $\models\varphi(a,b)$ and so 
$$(K',(\sigma'_g)_{g\in G})\models (\exists x)\,(\varphi(a,x)\,\wedge\,\bigwedge\limits_{g\in X}\sigma_g(x)=x),$$
where $X$ denotes the finite set of generators of $G$.
Hence
$$(K,(\sigma_g)_{g\in G})\models (\exists x)\,(\varphi(a,x)\,\wedge\,\bigwedge\limits_{g\in X}\sigma_g(x)=x)$$
and for some $b_0\in K^G$ we have that $\models\varphi(a,b_0)$.
\end{proof}

Therefore we see that an existentially closed substructure with $G$-action has a quite tame substructure of invariants. The next subsection is dedicated to the converse of this implication, so we would like to show that 
``if the 
substructure of invariants is tame
then the whole substructure with $G$-action 
is existentially closed''.

\begin{remark}
In Proposition 3.56 from \cite{Hoff3}, it was shown that if $(K,\bar{\sigma})$ is an existentially closed substructure with $G$-action, then $K$ is PAC.
However, the aforementioned proposition assumes quantifier elimination, elimination of imaginaries and stability (but $G$ there can be arbitrary).
\end{remark}

\subsection{Invariants of existentially closed actions}

\begin{lemma}\label{lemma:invariants.regular}
(QE, FS) Assume that $G$ is finite, $(K,(\sigma_g)_{g\in G})\subseteq (K',(\sigma'_g)_{g\in G})$ is an extension of substructures with $G$-action, the group action of $G$ on $K$ is faithful and $\dcl(K)=K$.
If $K^G$ is algebraically $K$-strongly PAC, then $K^G\subseteq (K')^G$ is regular.
\end{lemma}

\begin{proof}
If $\dcl(K)=K$ then also $\dcl(K^G)=K^G$.
Moreover, $K^G\subseteq (K')^G$ is regular if and only if $K^G\subseteq \dcl((K')^G)$ is regular and there is a unique way of extending $G$-action from $K'$ over $\dcl(K')$. Therefore, without loss of generality, we assume that $K'=\dcl(K')$ and so $\dcl((K')^G)=(K')^G$.
We need to show that $(K')^G\cap\acl(K^G)=K^G$.

Let $a\in (K')^G\cap\acl(K^G)\setminus K^G$.
Because for every $g\in G$, we have that $\sigma_g\big(\tp(a/K)\big)=\tp\big(\sigma_g(a)/K\big)$ and $a\in (K')^G$, we see that $\tp(a/K)$ is a $G$-invariant type.
By Fact \ref{fact:2.10} and Lemma \ref{lemma.transitively}, we see that $\tp(a/K)$ is a unique extension of $\tp(a/K^G)$ over $K$.

As $a\in\acl(K^G)$ and $\acl(K^G)\ind_{K^G}K$ (e.g. Remark 5.3 in \cite{casasimpl}),
$\tp(a/K^G)\subseteq\tp(a/K)$ is a non-forking extension.
Because $K^G$ is algebraically $K$-strongly PAC, $\tp(a/K^G)$ is finitely satisfiable in $K^G$.
As $a\in\acl(K^G)$, this means that it must be $a\in K^G$.
\end{proof}

\begin{definition}\label{def:old.3.18}
Assume that $C\subseteq K\subseteq \FC$ and that $G$ is finite. We call the pair $(C,K)$ \emph{$G$-closed} if
$C\subseteq K$ is a Galois extension, $G(K/C)\cong G$ and 
there is no $K'\subseteq \acl(K)$, $K\subsetneq K'$, such that the action of $G(K/C)$ extends over $K'$.
\end{definition}

\begin{lemma}\label{lemma:Gclosed.Kpac}
(QE, FS) Assume that $G$ is finite, $(K,\bar{\sigma})$ is a substructure with $G$-action such that
action of $G$ on $K$ is faithful and $\dcl(K)=K$.
Then $(K^G,K)$ is $G$-closed if and only if $K^G$ is algebraically $K$-strongly PAC.
\end{lemma}

\begin{proof}
By Fact \ref{fact:2.10}, $K\subseteq\acl(K^G)$, $K^G\subseteq K$ is Galois and $G(K/K^G)\cong G$.

Assume that $(K^G,K)$ is $G$-closed and let $p(x)\in S(K^G)$ be algebraic with a unique extension $\tilde{p}(x)$ over $K$ (being a non-forking extension follows naturally from $\acl(K^G)\ind_{K^G}K$, e.g. Remark 5.3 in \cite{casasimpl}).
We have that $\tilde{p}$ is $G$-invariant and so, by Lemma \ref{lemma:extending},
if $b\models\tilde{p}$ then
there exists an extension of substructures with a $G$-action,
$$(K,(\sigma_g)_{g\in G})\subseteq (K',(\sigma'_g)_{g\in G})$$
such that $K'=\dcl(K,b)$ and $b\in (K')^G$.
As $K'=\dcl(K,b)\subseteq \acl(K^G)=\acl(K)$, it must be that $K=K'$, so $b\in K$ and finally $b\in K^G$.

Now, we show the right-to-left implication.
Assume that $K'\subseteq \acl(K)$ and there is an extension of substructures with $G$-action:
$$(K,(\sigma_g)_{g\in G})\subseteq (K',(\sigma'_g)_{g\in G}).$$
By Lemma \ref{lemma:invariants.regular}, $K^G\subseteq (K')^G$ is regular.
As $(K')^G\subseteq K'\subseteq\acl(K)=\acl(K^G)$ it must be
$(K')^G\subseteq\dcl(K^G)=K^G$, so $K^G=(K')^G$.
By the proof of Proposition 4.1 from \cite{DHL1}
and the Galois correspondence for finite extensions (e.g. Theorem 12 in \cite{invitation}),
there exists a finite tuple $b$ from $K$ such that $K=\dcl(K^G,b)$.
Moreover, by the same proof of Proposition 4.1 from \cite{DHL1}, we also have that 
$K'=\dcl((K')^G,b)$.
Because $K^G=(K')^G$, we have that $K=\dcl(K^G,b)=\dcl((K')^G,b)=K'$.
\end{proof}

The following remark is not important for the main results of this paper and its purpose is mainly to generalize Theorem 3.25 from \cite{HK3}. As we use in its proof the Elementary Equivalence for PAC structures (\cite{DHL1}), we need to add more assumptions.

\begin{remark}
Let $T$ be stable with elimination of quantifiers and elimination of imaginaries.
Assume that PAC is a first order property.
Suppose that $(C,K)\subseteq(C',K')$ is an extension of $G$-closed substructures
such that $C$ and $C'$ are PAC.  
Then $C\preceq C'$.
\end{remark}

\begin{proof}
It is enough to reproduce the proof of Theorem 3.25 from \cite{HK3}, but in this more general context.
By the proof of Theorem 3.22 from \cite{HK3} or more similar Lemma 3.54 from \cite{Hoff3}, we have that
$C$ and $C'$ are bounded PAC structures.
Thus, by Corollary 3.11 from \cite{DHL1}, it is enough to show that the restriction map $r:G(C')\to G(C)$ is an isomorphism. After combining Lemma \ref{lemma:Gclosed.Kpac} and Lemma \ref{lemma:invariants.regular}, we obtain that $C\subseteq C'$ is regular, so $r$ is an epimorphism.

By Theorem 4.4 from \cite{Hoff4}, $G(C)$ is projective, which means that there exists
embedding $h$ as in the following diagram
$$\xymatrix{G(C) \ar[r]^{=} \ar@{-->}_{h}[dr] & G(C) \\
& G(C') \ar[u]^{r}}$$
But then $\mathcal{G}_0:=h[G(C)]\leqslant G(C')$ is a closed subgroup such that
$r|_{\mathcal{G}_0}:\mathcal{G}_0\to G(C)$ is an isomorphism. 

Because $K\subseteq\acl(C)$ and $K'\subseteq\acl(C')$, the restriction maps $G(C)\to G$ and $G(C')\to G$ lead to the following commutative diagram
$$\xymatrix{G(C')\ar[dr] \ar[rr]^{r}& & G(C) \ar[dl] \\
& G &}$$
and so $\mathcal{G}_0\mathcal{N}=G(C')$ for $\mathcal{N}:=\ker\big(G(C')\to G\big)$.
By Lemma 3.31 from \cite{Hoff3}, this implies that $\mathcal{G}_0=G(C')$ as expected.
\end{proof}

\begin{theorem}\label{thm.4conditions}
(QE, FS, ST$+$) Assume that $G$ is finite, say $|G|=l$.
Let $(K,\bar{\sigma})$ be a substructure with $G$-action such that $G$ acts faithfully on $K$ and $\dcl(K)=K$.
The following are equivalent:
\begin{enumerate}
    \item 
    $(K,\bar{\sigma})$ is existentially closed,
    
    \item
    $K^G$ is $K$-strongly PAC,
    
    \item
    $K^G$ is PAC and algebraically $K$-strongly PAC,
    
    \item
    $K^G$ is PAC and $(K^G,K)$ is $G$-closed.
\end{enumerate}
\end{theorem}

\begin{proof}
By Lemma \ref{lemma:Gclosed.Kpac}, (3)$\iff$(4).
(1)$\Rightarrow$(2) follows by Lemma \ref{lemma:Kpac1}.
The implication (2)$\Rightarrow$(3) follows by definitions.
To get the theorem, we will show that (3)$\Rightarrow$(1).

Assume that $\dcl(K)=K$, the group action is faithful and that $K^G$ is PAC and algebraically $K$-strongly PAC.
Using Fact \ref{fact:2.10}, we obtain the following
\begin{itemize}
    \item $K\subseteq\acl(K^G)$,
    \item $K^G\subseteq K$ is a Galois extension,
    \item $G(K/K^G)\cong G$.
\end{itemize}
The proof of Proposition 4.1 from \cite{DHL1} gives us existence of 
a finite tuple $\bar{b}=(b_0,\ldots,b_{l-1})$ from $K$ such that
$K=\dcl(K^G,\bar{b})$.

Consider $(K,(\sigma_g)_{g\in G})\subseteq (K',(\sigma'_g)_{g\in G})$.
Without loss of generality, we may assume that $(K',(\sigma'_g)_{g\in G})$ is existentially closed, in particular $\dcl(K')=K'$. We have that the group action of $G$ on $K'$ is faithful, thus by Fact \ref{fact:2.10},
we have that $K'\subseteq\acl((K')^G)$, $(K')^G\subseteq K'$ is Galois, 
and $G(K'/(K')^G)\cong G$.
Lemma \ref{lemma:invariants.regular} gives us that $K^G\subseteq (K')^G$ is regular, which means that the restriction map $G(K'/(K')^G)\to G(K/K^G)$ is onto, and so it is an isomorphism of finite groups.
The last thing implies $$K'=\dcl((K')^G,\bar{b}).$$

Let $\bar{B}$ be some enumeration of $\{\sigma_g(b_i)\;|\;g\in G,\;i<l\}$.
We have that $K'=\dcl((K')^G,\bar{b})=\dcl((K')^G,\bar{B})$.
Assume that
$$(K',(\sigma'_g)_{g\in G})\models \phi(a)$$
for some tuple $a$ from $K'$ and some quantifier-free formula $\phi(x)\in L_G(K)$.
First, we may present $\phi(a)$ as $\varphi_0(\sigma'_{g_0}(a),\ldots,\sigma'_{g_{l-1}}(a))$, 
where $\varphi_0(x_0,\ldots,x_{l-1})\in L(K)$ is quantifier-free.
Second, since $K=\dcl(K^G,\bar{B})$, we may present 
$\varphi_0(\sigma'_{g_0}(a),\ldots,\sigma'_{g_{l-1}}(a))$ as 
$\varphi(\sigma'_{g_0}(a),\ldots,\sigma'_{g_{l-1}}(a),\bar{B})$,
where $\varphi(x_0,\ldots,x_{l-1},\bar{y})\in L(K^G)$ is quantifier-free.

Let $\sigma'_{g_0}=\id_L$, so $\sigma'_{g_0}(a)=a$.
Because $a\in K'=\dcl((K')^G,\bar{B})$, there exists a finite tuple $\bar{c}\subseteq (K')^G$ and a quantifier-free formula $\psi_0(\bar{z},\bar{y},x)\in L$ such that
\begin{itemize}
    \item $\psi_0(\bar{c},\bar{B},\FC)=\{a\}$,
    \item $\models (\forall\bar{z},\bar{y},x,x')\,\big(\psi_0(\bar{z},\bar{y},x)\,\wedge\,\psi_0(\bar{z},\bar{y},x')\,\longrightarrow\,x=x'\big)$.
\end{itemize}
Because $\sigma_{g_i}$ permutes $\bar{B}$, there exists a permutation $s_i$ such that
$\sigma_{g_i}(\bar{B})=s_i(\bar{B})$. We define $\psi_i(\bar{z},\bar{y},x)$ as $\psi_0(\bar{z},s_i(\bar{y}),x)$.
Note that $\psi_i(\bar{c},\bar{B},\FC)=\{\sigma'_{g_i}(a)\}$ and
$$(K',(\sigma_g')_{g\in G})\models $$
$$(\forall\bar{z},x,x')\,\big(\bigwedge\limits_{g\in G}\sigma_g(\bar{z})=\bar{z}\,\wedge\,\psi_0(\bar{z},\bar{B},x)\,\wedge\,\psi_i(\bar{z},\bar{B},x')\,\rightarrow\,\sigma_{g_i}(x)=x'\big).$$
To see the last line, let $\bar{d}\subseteq (K')^G$, $m,m'\in K'$ be such that
$$\models \psi_0(\bar{d},\bar{B},m)\,\wedge\,\psi_i(\bar{d},\bar{B},m').$$
We do know that $\psi_0(\bar{d},\bar{B},\FC)=\{m\}$, which after applying an extension $\tilde{\sigma}_{g_i}\in\aut(\FC)$ of $\sigma'_{g_i}$ changes it into $\psi_0(\bar{d},s_i(\bar{B}),\FC)=\{\sigma'_{g_i}(m)\}$.
We have that
$$m'\in\psi_i(\bar{d},\bar{B},\FC)=\{\sigma'_{g_i}(m)\}.$$

Since the whole formula is universal and has only parameters from $K$, it follows that
$$(K,(\sigma_g)_{g\in G})\models 
(\forall\bar{z},x,x')\,\big(\bigwedge\limits_{g\in G}\sigma_g(\bar{z})=\bar{z}\,\wedge\,\psi_0(\bar{z},\bar{B},x)\,\wedge\,\psi_i(\bar{z},\bar{B},x')\,\rightarrow\,\sigma_{g_i}(x)=x'\big),$$
where $i<l$.

Consider $p(\bar{z}):=\tp(\bar{c}/K^G)$.
Because $K^G\subseteq (K')^G$ is regular (thus also primary) and $\bar{c}\subseteq (K')^G$, 
Fact \ref{fact.3.35} implies that $p(\bar{z})$ 
is stationary.
As $K^G$ is PAC, the type $p(\bar{z})$ is finitely satisfiable in $K^G$.
The tuple $\bar{B}\subseteq K$ is algebraic over $K^G$, hence there exists a quantifier-free $\theta(\bar{y})\in L(K^G)$ such that $\theta(\bar{y})\vdash\tp(\bar{B}/K^G)$.
The following formula
$$(\exists\,\bar{y},\,x_0,\ldots,x_{l-1})\,\big( \bigwedge\limits_{i<l}\psi_i(\bar{z},\bar{y},x_i)\,
\wedge\,\varphi(x_0,\ldots,x_{l-1},\bar{y})\,
\wedge\,\theta(\bar{y})\big)$$
belongs to $p(\bar{z})$, thus there exists $\bar{d}\subseteq K^G$ such that
$$\models(\exists\,\bar{y},\,x_0,\ldots,x_{l-1})\,\big( \bigwedge\limits_{l<e}\psi_i(\bar{d},\bar{y},x_i)\,
\wedge\,\varphi(x_0,\ldots,x_{l-1},\bar{y})\,
\wedge\,\theta(\bar{y})\big).$$
It means that there are $\bar{B}'\subseteq\FC$ and $a'_0,\ldots,a'_{l-1}\in \FC$ such that
$$\models  \bigwedge\limits_{i<l}\psi_i(\bar{d},\bar{B}',a'_i)\,
\wedge\,\varphi(a'_0,\ldots,a'_{l-1},\bar{B}')\,
\wedge\,\theta(\bar{B}').$$
Since $\models \theta(\bar{B})\,\wedge\,\theta(\bar{B}')$, there exists $f\in\aut(\FC/K^G)$ such that $f(\bar{B}')=\bar{B}$. By applying $f$, we obtain
$$\models  \bigwedge\limits_{i<l}\psi_i(\bar{d},\bar{B},f(a'_i))\,
\wedge\,\varphi(f(a'_0),\ldots,f(a'_{l-1}),\bar{B}).$$
We have that $\psi_0(\bar{d},\bar{B},\FC)=\{f(a'_0)\}$.
Since, for each $i<l$, the subset $\psi_i(\bar{d},\bar{B},\FC)=\psi_0(\bar{d},s_i(\bar{B}),\FC)$ is an automorphic image of $\psi_0(\bar{d},\bar{B},\FC)$, it must be that $|\psi_i(\bar{d},\bar{B},\FC)|=1$ and so $f(a'_i)\in \dcl(K^G,\bar{B})=K$ for each $i<l$.
Moreover, we have that $\sigma_{g_i}(f(a'_0))=f(a'_i)$ for each $i<l$.
Therefore $\models \varphi(f(a'_0),\ldots,f(a'_{l-1}),\bar{B})$ leads to
$$(K,(\sigma_g)_{g\in G})\models \varphi\big(\sigma_{g_0}(f(a'_0)),\ldots,\sigma_{g_{l-1}}(f(a'_{0})),\bar{B}\big),$$
as expected.
\end{proof}

\begin{cor}
(QE, FS, ST$+$) Let $G$ be finite and let $(K,\bar{\sigma})$ be a substructure with $G$-action. Then $(K,\bar{\sigma})$ is existentially closed if and only if
\begin{enumerate}
    \item $\dcl(K)=K$,
    \item the group action of $G$ on $K$ is faithful and,
    \item $K^G$ is $K$-strongly PAC.
\end{enumerate}
\end{cor}

\begin{proof}
If the conditions (1), (2) and (3) hold
then 
$(K,\bar{\sigma})$ is existentially closed by Theorem \ref{thm.4conditions}.

If $(K,\bar{\sigma})$ is existentially closed, 
then 
Remark \ref{rem:dcl} gives us that $\dcl(K)=K$, 
Lemma \ref{lemma:faithful} gives that the group action is faithful,
and the fact that $K^G$ is $K$-strongly PAC follows from Lemma \ref{lemma:Kpac1}.
\end{proof}

In a very similar way, we conclude the following.

\begin{cor}\label{cor:ec}
(QE, FS, ST$+$) Let $G$ be finite and let $(K,\bar{\sigma})$ be a substructure with $G$-action. Then $(K,\bar{\sigma})$ is existentially closed if and only if
\begin{enumerate}
    \item $\dcl(K)=K$,
    \item the group action of $G$ on $K$ is faithful and,
    \item $K^G$ is PAC and algebraically $K$-strongly PAC.
\end{enumerate}
\end{cor}

\subsection{Existence of model companion}

\begin{remark}\label{rem:unique.vs.invariant}
Assume that $A\subseteq C$ is a Galois extension (e.g.:
if QE, FS, $G$ is finite and $G$ acts faithfully on $K=\dcl(K)$, we can take $A=K^G$ and $C=K$).
Let $p(x)\in S(A)$. The following are equivalent:
\begin{enumerate}
    \item There exists unique extension of $p$ over $C$.
    \item There exists a $G(C/A)$-invariant extension of $p$ over $C$.
\end{enumerate}
\end{remark}

\begin{proof}
If there is only one extension of $p(x)$ over $C$ it is automatically $G(C/A)$-invariant.
Assume that $p(x)$ has a $G(C/A)$-invariant extension over $C$, say $p_1(x)\in S(C)$
and led $p_2(x)\in S(C)$ be also an extension of $p(x)$. 
As $p_1|_{A}=p_2|_{A}$, there exists $f\in\aut(\FC/A)$ such that $f(p_1)=p_2$.
Since $A\subseteq C$ is Galois, we know that $f|_C=\sigma$ for some $\sigma\in G(C/A)$.
Then, $p_2=f(p_1)=f|_C(p_1)=\sigma(p_1)=p_1$.
\end{proof}

\begin{definition}
Let $\pi(x)$ be a partial type over $A$.
We say that $\pi$ is \emph{$A$-irreducible} if there exists $p(x)\in S(A)$ such that $\pi\vdash p$.
\end{definition}

\begin{remark}\label{rem:irreducibility}
Let $\varphi(x)\in L(A)$ be a consistent formula.
Then, $\varphi(x)$ is \emph{$A$-irreducible}
if and only if $\{\varphi(x)\}$ is $A$-irreducible, which is equivalent to saying that
there are no formulae $\varphi_1(x),\varphi_2(x)\in L(A)$ such that $\varphi(\FC)\,\cap\,\varphi_1(\FC)\neq\emptyset$, $\varphi(\FC)\,\cap\,\varphi_2(\FC)\neq\emptyset$ and $\varphi(\FC)\,\cap\,\varphi_1(\FC)\,\cap\,\varphi_2(\FC)=\emptyset$. We use this characterization in the crucial Remark \ref{rem:FO.Kirred}.
\end{remark}

\begin{lemma}\label{lemma:algKpac2}
(QE, FS) Assume that $G$ is finite and $G$ acts faithfully on $K=\dcl(K)$.
The following are equivalent.
\begin{enumerate}
    \item 
    $K^G$ is algebraically $K$-strongly PAC.
    
    \item 
    Each algebraic type $p(x)\in S(K^G)$, which has a $G$-invariant extension over $K$, is satisfiable in $K^G$.
    
    \item
    Each $G$-invariant algebraic type $\tilde{p}(x)\in S(K)$ is satisfiable in $K$.
    
    \item
    For each $\theta(x)\in L(K^G)$, if $0<|\theta(\FC)|<\omega$ and $\theta(\FC)$ is $K$-irreducible then $\theta(K^G)\neq\emptyset$.
\end{enumerate}
\end{lemma}

\begin{proof}
The proof is easy, so we only sketch it.
By Remark \ref{rem:unique.vs.invariant}, we immediately obtain (1)$\iff$(2).
For the (2)$\iff$(3), it is enough to observe that a $G$-invariant algebraic type is isolated by a $L(K^G)$-formula,
so its restriction to $K^G$ is also algebraic.
We argue similarly on (3)$\iff$(4): a $G$-invariant algebraic type over $K$ is isolated by an $L(K^G)$-formula,
which is consistent, algebraic (i.e. has finitely many realizations) and $K$-irreducible.
\end{proof}

\begin{remark}\label{rem:FO.Kirred}
(QE, FS$+$) In this remark, we investigate in what way being a $K$-irreducible formula may be expressed as a first order statement.
Assume that $K=\dcl(K)$ and consider a quantifier-free formula $\varphi(y,x)\in L$ and a tuple $a\in K^y$.
Moreover assume that $0<|\varphi(a,\FC)|=n<\omega$.
Recall that we can use the technical condition (FS$+$), listed at the very beginning of the paper.
Let $S_x$ be the sort related to the variable $x$ and let $E(\bar{x},\bar{x}')$ be a $\emptyset$-definable equivalence relation given by a formula expressing that ``$\{x_1,\ldots,x_{n-1}\}=\{x'_1,\ldots,x'_{n-1}\}$'', where $\bar{x}=(x_1,\ldots,x_{n-1})$ and $\bar{x}'=(x'_1,\ldots,x'_{n-1})$ are tuples of variables from $S_x$. 
As we assume, $E$ is the fibration of a $0$-definable function $f:(S_x)^{n-1}\to S_w$.
Note that the elements of the image of $f$ correspond to nonempty subsets of $S_x(\FC)$ of the size at most $n-1$.

One more thing before coming to the point. Assume that the element $d$ belongs to the sort $S_w(\FC)$,
then the formula
$$(\exists\, x_1,\ldots,x_{n-1}\in S_x)\,\big( f(x_1,\ldots,x_{n-1})=d\,\wedge\,\bigwedge\limits_{i=1}^{n-1}\varphi(a,x_i)\big)$$
is modulo $T$ equivalent to a quantifier free formula, say $\xi_{\varphi,n}(a,d)$.

Now, $\varphi(a,x)$ is $K$-irreducible if and only if there is no proper subset $\emptyset\neq X\subsetneq\varphi(a,\FC)$, such that $X$ is $K$-definable. 
In other words, for each proper subset $\emptyset\neq X\subsetneq\varphi(a,\FC)$, we have that the code $\ulcorner X\urcorner$ does not belong to $\dcl(K)=K$.
We can express this last sentence, in the structure $K$, as follows
$$K\models\neg(\exists\,w\in S_w)\,(\xi_{\varphi,n}(a,w)).$$
We will use the above in the proof of Theorem \ref{thm:ec.if.PAC}.
\end{remark}

Now, we want to show that being a definably closed (in $L$) subset of $\FC$ is a first order statement.

\begin{remark}\label{rem:QE}
(QE) Consider any $\theta(y,x)\in L$. The formula
$$\theta(y,x)\,\wedge\,(\forall\, x_1,x_2)\,\big(\theta(y,x_1)\,\wedge\,\theta(y,x_2)\,\rightarrow\,x_1=x_2\big)$$
is equivalent modulo $T$ to some quantifier-free $\psi_{\theta}(y,x)\in L$. Moreover, also the formula $(\exists\,x)\,(\psi_{\theta}(y,x))$ is equivalent modulo $T$ to a quantifier-free formula $\psi^0_{\theta}(y)\in L$.
Consider 
$$\Sigma:=\{\psi^0_{\theta}(y)\,\rightarrow\,(\exists\,x)\,(\psi_{\theta}(y,x)\big)\;|\;\theta(y,x)\in L\}.$$
\end{remark}

\begin{lemma}\label{lemma:dcl}
(QE, in the notation of Remark \ref{rem:QE}) For a substructure $B\subseteq\FC$,
$B\models\Sigma$ if and only if $\dcl(B)=B$.
\end{lemma}

\begin{proof}
Assume that $B\models\Sigma$ and let $b\in\dcl(B)$.
There exists a formula $\theta(y,x)\in L$ and $a\in B$ such that $\theta(a,\FC)=\{b\}$.
Then $\models \psi_{\theta}(a,b)$ and $\models\psi^0_{\theta}(a)$.
As $\psi^0_{\theta}$ is quantifier-free, also $B\models \psi^0_{\theta}(a)$, thus $B\models (\exists\,x)\,(\psi_{\theta}(a,x))$. It means that there exists $b'\in B$ such that $\models\psi_{\theta}(a,b')$.
We see that $b'=b$ and so $b\in B$.

Now, let $B=\dcl(B)$. Assume that $B\models \psi^0_{\theta}(a)$ for some $a\in B$ and $\theta(y,x)\in L$.
We have $\models\psi^0_{\theta}(a)$, so there exists some $b\in \FC$ such that $\models\psi_{\theta}(a,b)$.
This implies that $b\in\dcl(a)\subseteq \dcl(B)=B$. Therefore there exists $b\in B$ such that $B\models\psi_{\theta}(a,b)$.
\end{proof}

The following theorem is an answer towards Question 2.9 (also Question 5.1) and Conjecture 5.2 from \cite{Hoff3}.

\begin{theorem}\label{thm:ec.if.PAC}
(QE, FS$+$, ST$+$) Let $G$ be finite.
The model companion of the theory of substructures with $G$-action exists provided
PAC is a first order property.
\end{theorem}

\begin{proof}
By Corollary \ref{cor:ec} and Lemma \ref{lemma:algKpac2}, 
we need to write down as first order statements the following conditions:
\begin{enumerate}
    \item[(0)] $(K,(\sigma_g)_{g\in G})$ is a substructure with $G$-action,
    \item $\dcl(K)=K$,
    \item the group action of $G$ on $K$ is faithful,
    \item $K^G$ is PAC,
    \item for each $\theta(x)\in L(K^G)$, if $0<|\theta(\FC)|<\omega$ and $\theta(\FC)$ is $K$-irreducible then $\theta(K^G)\neq\emptyset$.
\end{enumerate}
We are working in the language $L_G$.
The condition (0) is naturally a first order statement, similarly the condition (2).
Lemma \ref{lemma:dcl} shows that also the condition (1) is a first order statement.
By the assumptions, condition (3) is a first order statement.
To finish the proof of the theorem,
we need to show that the condition (4) is also a first order statement.

There is no harm in assuming that the formula $\theta(x)$ is $\varphi(a,x)$ for some tuple $a$ from $K^G$ and some quantifier-free formula $\varphi(y,x)\in L$.
The condition (4) will be expressed as an axiom scheme running over all quantifier-free formulae $\varphi(y,x)\in L$ and all $0<n<\omega$.

Fix a quantifier-free formula $\varphi(y,x)\in L$ and a natural number $n>0$.
There exists a quantifier-free $L$-formula $\psi_{\varphi}(y)$ equivalent modulo $T$ to the formula $(\exists^{=n}\,x)\,(\varphi(y,x))$. 
We are in situation of Remark \ref{rem:FO.Kirred}, so we can involve the formula $\xi_{\varphi,n}(y,w)$.
Our axiom scheme may be written as:
$$
(\forall\,y)\,\Big(\bigwedge\limits_{g\in G}\sigma_g(y)=y\,\wedge\,\psi_{\varphi}(y)\,\wedge\,
\neg(\exists\,w\in (S_x)^{n-1}/E)\,(\xi_{\varphi,n}(y,w))$$
$$\rightarrow\;(\exists\,x)\big(\bigwedge\limits_{g\in G}\sigma_g(x)=x\,\wedge\,\varphi(y,x)\big)\Big).$$
\end{proof}

\begin{question}\label{Q324}
(QE, FS$+$, ST$+$) 
Can we obtain a converse of Theorem \ref{thm:ec.if.PAC}?
More precisely, does the following equivalence hold:
the model companion of the theory of substructures with $G$-action exists for every finite group $G$ 
if and only if PAC is a first order property?
\end{question}

\begin{remark}	
After writing the proofs of Theorems \ref{thm.4conditions} and \ref{thm:ec.if.PAC}, we have noticed (but we have not checked all the details) that this result holds in a much greater generality, that is: if in the definition of PAC we replace ``stationary'' with ``acl-stationary'' (a unique extension over algebraic closure of the parameters), then the assumptions of stability, coding finite sets, and eliminating strong types may be skipped in Theorem \ref{thm:ec.if.PAC}. However, in this case it is unclear how useful such a result would be in terms of axiomatizing existentially closed finite group actions in this case, since there is no guarantee that faithful actions of a finite groups exist at all in general (consider for example the theory of linear orders) and the faithfulness in the stable was guaranteed by Lemma \ref{lemma:faithful}.
\end{remark}

\section{PAC structures in particular theories}\label{sec:examples}
In this section, we discuss the PAC property in some specific cases as well as some general methods for understanding PAC structures with respect to a given theory. As we are going to consider the notions of a regular extension and of a PAC structure in different theories, we plan to write ``$T$-regular'' and ``$T$-PAC'' instead of ``regular in $T$'' and ``PAC in $T$'' respectively.

We will often refer to several particular stable theories as: the theory of compact complex manifolds CCM (for background, the reader is referred to \cite{Moosacpx}) and the theories of differentially closed fields of characteristic $0$ denoted $\mathrm{DCF}_0$ (see e.g. \cite{MMPfieldsMarker}) and its positive characteristic version $\mathrm{DCF}_p$ (see e.g. \cite{Wood1} and \cite{Wood2}),
and the theories of separably closed fields of positive characteristic $\mathrm{SCF}_{p,e}$ and
$\mathrm{SCF}_{p,\infty}$ (see e.g. \cite{messmer}).

\subsection{General methods}\label{secgen}
In this subsection, we focus on two general contexts in which the PAC property is well understood. However, in both these cases showing that PAC is a first-order property requires some extra work.

\subsubsection{Totally transcendental theories}
In this part, we assume that the theory $T$ is $\omega$-stable. As before, let us fix for convenience a monster model $\mathfrak{C}$ of $T$ and an arbitrary small substructure $K\subset \mathfrak{C}$. It is well-known that stationary types in $\omega$-stable theories are determined by the formulas of Morley degree one belonging to them. In particular, we have the following result, which actually coincides with Hrushovski's definition of the PAC property in the strongly minimal case (see \cite[Definition 1.2]{manuscript} and \cite[Proposition 3.10]{Hoff3}).
\begin{prop}\label{dmp}
If $T$ is a $\omega$-stable theory, then $K$ is $T$-PAC if and only if for any formula $\varphi\in L(K)$ of multiplicity (Morley degree) one, we have that $\varphi(K)\neq \emptyset$.
\end{prop}
We recall that ``DMP'' stands for ``Definable Multiplicity Property'' and it says that for any formula $\phi(x;a)\in L(K)$, there is a formula $\theta(y)\in \tp(a)$ such that whenever
$\mathfrak{C}\models \theta(a')$ then we have:
$$\mathrm{RM}\left(\phi(x;a')\right)=\mathrm{RM}\left(\phi(x;a)\right),\ \ \ \ \deg_M\left(\phi(x;a')\right)=\deg_M\left(\phi(x;a)\right)$$
  (see e.g. \cite[Definition 1.1]{KiPi}). Some $\omega$-stable theories have DMP and some do not (see Remark \ref{dmprem} below).
We get the following obvious conclusion, which was also stated in \cite{Afshordel} under the assumption of finiteness of the Morley rank.
\begin{prop}\label{dmppac}
If $T$ is $\omega$-stable with quantifier elimination and has DMP, then being $T$-PAC is first-order.
\end{prop}
\begin{proof}
Since $T$ has DMP, for each $\phi(x;y)\in L$, there is $\theta_{\phi}(y)$ such that for all $c\in \mathfrak{C}^{|y|}$, we have:
$$\mathfrak{C}\models \theta_{\phi}(c)\ \ \ \ \ \text{if and only if}\ \ \ \ \ \deg_M\left(\phi(x;c)\right)=1.$$
Therefore, it is easy to write down a first-order axiom scheme expressing the $T$-PAC property.
\end{proof}
\begin{remark}\label{dmprem}
We comment here on several particular $\omega$-stable theories.
\begin{enumerate}
\item Proposition \ref{dmppac} applies to the case of $T=\mathrm{ACF}_p$, that is to the classical notion of PAC.

\item It is known that Morley degree is not definable in the theory DCF$_0$ (see \cite[Question 1.2]{Freitag1} and \cite{FSLi}). However, DCF$_0$-PAC is still first-order as it was shown in \cite{PilPol}.

\item It is open whether the theory of compact complex manifolds has DMP, however another approach towards the PAC property works here, which will be discussed in the next part. Partial results towards DMP for the theory CCM were obtained in \cite{Dale2}.

\end{enumerate}
\end{remark}

\subsubsection{Noetherian theories}\label{secnoeth}
In this part, we assume that models of $T$ are naturally equipped with an extra topological structure. This assumptions is modelled on the case of $T=\mathrm{ACF}_p$ and the Zariski topology. Such issues were thoroughly discussed in \cite{zilber10}. We diverge here a bit from the set-up of \cite{zilber10} to cover the case of the theory DCF$_0$ as well.

We start from a purely topological context. Assume that $S$ is a Noetherian topological space and let $\mathbb{B}$ be the Boolean algebra of constructible sets in $S$. In this part, ``irreducible'' always refers to the topological irreducibility with respect to a given Noetherian topology.
The following properties are folklore and they can be easily checked.
\begin{itemize}
  \item If $V$ is a non-empty closed irreducible subset of $S$, then
  $$p_V:=\{C\in \mathbb{B}\ |\ \mathrm{int}_V(C\cap V)\neq \emptyset\}$$
  is an ultrafilter on $\mathbb{B}$.
  \item The map $V\mapsto p_V$ is a bijection between the set of closed irreducible subsets of $S$ and the set of ultrafilters on $\mathbb{B}$.
\end{itemize}
We specify now our model-theoretic context.
\begin{definition}
By a \emph{Noetherian theory}, we mean a pair $(T,\sum)$, where $T$ is a complete $L$-theory and $\sum$ consists of $L$-formulas of the form $\varphi(x;y)$, where the variables $x,y$ vary, such that for any $M\models T$ and any $A\subseteq M$, we have the following.
\begin{itemize}
  \item A subset $V\subseteq M^{|x|}$ is said to be $A$-\emph{closed} if and only if there is $a\subset A$ and $\varphi(x;y)\in \sum$ such that $V=\varphi(M;a)$.
  \item The family of $A$-closed sets constitutes the family of closed sets of a Noetherian topology, which we call the \emph{$A$-topology}.
  \item Constructible sets with respect to the $A$-topology coincide with $A$-definable subsets (in Cartesian powers of $M$).
\end{itemize}
\end{definition}
\begin{remark}
\begin{enumerate}
    \item It should automatically follow (possibly after adding some light assumptions such as the equality being in $\sum$)
that models of our Noetherian theories are \emph{topological structures} in the sense of \cite[Definition 5.1]{BGH} and \cite[Section 2]{zilber10}.

    \item The referee has pointed out to us that a very similar notion of a Noetherian theory was recently introduced by Martin-Pizarro and Ziegler (see \cite[Definition 2.18]{MPZ}).
\end{enumerate}
\end{remark}

\begin{example}
We discuss several examples and non-examples of the above situation.
\begin{enumerate}
  \item The theory of algebraically closed fields (of a given characteristic) is Noetherian by considering the Zariski topology.

  \item  The theory of compact complex manifolds (CCM) is also Noetherian, where the (Zariski) Noetherian topology is given by closed analytic subsets (see \cite[Section 3.4.2]{zilber10}).


  \item In the case of differential fields, we have the Kolchin topology.
  \begin{itemize}
  \item The theory $\mathrm{DCF}_0$ is Noetherian by \cite[Theorem 2.4]{MMPfieldsMarker}.

  \item More generally, the theory $\mathrm{DCF}_{0,m}$ is Noetherian by \cite[Theorem 3.1.7]{McGrail}.



\item The theory $\mathrm{DCF}_p$ is not an example, since there is no quantifier elimination down to Kolchin constructible sets (see \cite[Section 3]{Wood1}).
  \end{itemize}

\item  The theory $\mathrm{SCF}_{p,e}$ with the $\lambda$-topology is not an example, since the $\lambda$-topology is not Noetherian (see \cite[Section 4.6]{messmer}).

\end{enumerate}
\end{example}
For a fixed $A\models T_{\forall}$ and $n>0$, it is clear that the map $V\mapsto p_V$ is a bijection between the set of appropriate $A$-closed $A$-irreducible sets and the Stone space $S_n(A)$ of $n$-types over $A$. In particular, any Noetherian theory is $\omega$-stable. We still need to have a connection between the topology and forking, which is given by the following.
\begin{prop}\label{topfork}
Assume that $A\subseteq M\models T$ and $p_V\in S_n(M)$. Then, the type $p_V$ does not fork over $A$ if and only if $V$ is definable over $\acl(A)$.
\end{prop}
\begin{proof}
Since $p_V$ does not fork over $A$ if and only if it does not fork over $\acl(A)$, we can and will assume that $A=\acl(A)$.
\\
$(\Rightarrow)$
Let $V=V_b$ and assume that $V_b$ is not definable over $A$. Let us define
$$V_0:=\bigcap_{\mathrm{tp}(c/A)=\mathrm{tp}(b/A)}V_c.$$
Since $V_b$ is not definable over $A$, we get that $V_0\subsetneq V$. By Noetherianity, $V$ is definable and closed. Since $V_0$ is $A$-invariant, we get that $V_0$ is $A$-definable. In particular, the formula ``$x\in V\setminus V_0$'' belongs to $p_V$. Since the formula ``$x\in V\setminus V_0$'' forks over $A$ (see e.g. the characterization of forking from \cite[Lemma 2.16(c)]{anandgeometric}), the type $p_V$ forks over $A$.
\\
$(\Leftarrow)$
 We assume that $V$ is $A$-definable. It is enough to show that for any proper $M$-closed $W=W_b\subset V$, we have that the formula ``$x\in V\setminus W$' does not fork over $A$. If this formula forks over $A$, then by (the logical) compactness, there is a finite set of $A$-conjugates $b=b_1,\ldots,b_n$ such that:
$$\left(V\setminus W_{b_1}\right)\cap \ldots \cap \left(V\setminus W_{b_n}\right)=\emptyset.$$
But then $V=W_{b_1}\cup \ldots \cup W_{b_n}$ and each $W_{b_i}$ is a proper $M$-closed subset of $V$, which contradicts the $M$-irreducibility of $V$.
\end{proof}
We obtain the expected description of stationary types.
\begin{cor}\label{irrstat}
Let $A\models T_{\forall}$,  and $p_W\in  S_n(A)$. Then, $p_W$ is stationary if and only if $W$ is \emph{absolutely irreducible}, that is for any $M\models T$ containing $A$ as a substructure, $W$ is irreducible in the $M$-topology.
\end{cor}
\begin{proof}
Let $W=W_1\cup \ldots \cup W_n$ be the decomposition of $W$ into the $M$-irreducible $M$-closed components. By uniqueness, each $W_i$ is defined over $\acl(A)$. By Proposition \ref{topfork}, each type $p_{W_i}$ does not fork over $A$. Since for each $i$, we have $\mathrm{cl}_A(W_i)=W$, we get that each $p_{W_i}$ extends $p_W$. It is easy to see now that $W_i$'s correspond exactly to non-forking extensions of $p_W$, which concludes the proof.
\end{proof}

Similarly as in the case of Proposition \ref{dmp}, we get the following result.
\begin{prop}\label{pacirr}
For any $K\models T_{\forall}$, we have that $K$ is $T$-PAC if and only if for any absolutely irreducible $K$-closed set $V$ and any non-empty relatively $K$-open $U\subseteq V$, we have that $U(K)\neq \emptyset$.
\end{prop}
\begin{remark}
\begin{enumerate}
\item In the cases of $T=\mathrm{ACF}_p$ and $T=\mathrm{DCF}_{0,m}$, we can just consider the condition ``$V(K)\neq \emptyset$'' in Proposition \ref{pacirr}, since these topologies have basis of open sets being definably isomorphic to affine closed sets.
It looks like there is no similar simplification for the theory CCM, since (at least in the category of complex manifolds) being isomorphic to a compact complex manifold would imply being closed.

\item  Proposition \ref{pacirr} together with Item $(1)$ above directly generalizes the classical case of $T=\mathrm{ACF}_p$.

\item For $T=\mathrm{DCF}_{0,m}$ the description from Proposition \ref{pacirr} (together with Item $(1)$ above) coincides with the \emph{definition} taken in \cite[Section 5.16]{SanchezTressl} (see also \cite[Remark 4.7(1)]{HoffSanch}).

\item For $T=\mathrm{CCM}$, we believe that this notion has not been considered before.
\end{enumerate}
\end{remark}
Similarly as in the previous part, we get the following result.
\begin{prop}\label{pacirrfo}
Assume that $T$ is a Noetherian theory. If the topological irreducibility is definable in $T$, then $T$-PAC is first-order.
\end{prop}
We would like to single out one important case below.

\begin{cor}\label{cor.PAC.CCM}
$\mathrm{CCM}$-PAC is first-order.
\end{cor}
\begin{proof}
By Proposition \ref{pacirrfo}, it is enough to show that the topological irreducibility is definable in the theory CCM, which was shown in \cite{Dale1}.
\end{proof}
\begin{remark}
We would like to mention that Rahim Moosa pointed out to us that an argument for the definability of the topological irreducibility in the case of compact complex spaces can be also found in an earlier work of Campana. The reader is advised to consult $\S$.3.B of Premiere Partie of \cite{Campanacycle}.
\end{remark}

For the terminology used in the next theorem, we advise the reader to recall Definition 3.1(5).

\begin{theorem}\label{thm:GCCM}
If $G$ is finite, then the theory $G$-$\mathrm{CCM}$, exists and it is supersimple with geometric elimination of imaginaries, codes finite sets and has ``semi" quantifier elimination (in the same way as the theory $\mathrm{ACFA}$).
\end{theorem}

\begin{proof}
The existence of $G$-CCM follows by Theorem \ref{thm:ec.if.PAC} and Corollary \ref{cor.PAC.CCM}.
The properties of $G$-CCM, listed in the statement, follow by Corollary 4.28, Theorem 4.36,
Lemma 4.37 and Remark 4.13 from \cite{Hoff3} after noticing that CCM is superstable with elimination of quantifiers and elimination of imaginaries.
\end{proof}

\begin{remark}
All the properties stated in Theorem \ref{thm:GCCM} hold also in the case of the theory CCMA (\cite{hils0}). One difference between $G$-CCM, for finite $G$, and CCMA are the values of the SU-rank. As ACFA is stably embedded in CCMA, there is a sort in CCMA on which the SU-rank is not finite. On the other hand, one can show that the SU-rank of $G$-CCM is finite (sort-by-sort).
\end{remark}

\begin{remark}\label{noethlast}
We finish this part with some comments on the theories of differential fields.
\begin{enumerate}

  \item Let us recall here that for the theory DCF$_{0,m}$ the definability of Kolchin topological irreducibility is equivalent to the (generalized) \emph{Ritt problem} (see \cite{FSLi}). However, the PAC property for DCF$_{0,m}$ is still first-order, which was shown in \cite{SanchezTressl}.

  \item It may be a good moment to point out that the methods of this paper cover that $\mathrm{DCF}_{0}$-PAC is first-order, but fail to generalize it to the case of DCF$_{0,m}$ for $m>1$. In \cite{SanchezTressl}, the general case is shown in the following way.
      \begin{itemize}
        \item First the authors of \cite{SanchezTressl} show that ``differential largeness'' is a first-order property (see \cite[Proposition 4.7]{SanchezTressl}).

        \item Then they show that DCF$_{0,m}$-PAC is equivalent with the classical field PAC together with the ``differential largeness'' (see \cite[Section 5.16]{SanchezTressl}).
      \end{itemize}
      The above scheme of a proof looks like a possible another general approach (at least for the theories of fields with operators).
      We will discuss it further at the end of this paper (see Remark \ref{last}).
\end{enumerate}
\end{remark}

\subsubsection{Equational theories}\label{secequa}

In \cite{PiSr}, the notion of an \emph{equational theory} is introduced. Briefly, a theory $T$ is equational if any formula is equivalent (modulo $T$) to a Boolean combination of instances of equations, where a formula $\varphi(x,y)$ is an equation (modulo $T$) if the family of definable sets given by finite intersections of its instances (in any model of $T$) has the DCC (Descending Chain Condition).

The set-up of equation theories generalizes, in some sense, the set-up of  Noetherian theories from Section \ref{secnoeth}, since Noetherian theories are equational ``in a strong sense'' that is the DCC condition holds not only in the case of instances of one formula but for all closed sets with respect to the given Noetherian topology. Not all the equational theories are Noetherian, since any Noetherian theory is $\omega$-stable and, for example, $\mathrm{Th}(\Zz,+)$ is equational and not $\omega$-stable (see Remark at the end of Section 2 in \cite{PiSr}).

There is a natural notion of irreducibility in the case of equational theories and it is possible that for an equation theory $T$, if this notion of irreducibility is definable, then $T$-PAC is first order.

\subsection{Fields of positive characteristic}
In this subsection, we focus on three stable theories of fields of positive characteristic: $\mathrm{SCF}_{p,e}$ ($e$ finite), $\mathrm{SCF}_{p,\infty}$, and $\mathrm{DCF}_p$.
There are several possible languages to consider for these fields of positive characteristic and we will actually use three different options here.

\subsubsection{Separably closed fields of finite imperfection degree}
We consider the theory $\mathrm{SCF}_{p,e}$ ($e$ finite) in the language $L_{\lambda,b}$, where $b$ stands for an $e$-tuple of constant symbols corresponding to a fixed $p$-basis and we also have symbols for unary $\lambda$-functions defined with respect to $b$ (see \cite[Section 1.8]{Cha02}). Then, the theory SCF$_{p,e}$ has quantifier elimination and elimination of imaginaries (cf. \cite[Section 1.8]{Cha02}).

Let us fix as usual a monster model $\mathfrak{C}\models \mathrm{SCF}_{p,e}$. We make the following identification $b=b^{\mathfrak{C}}$. A subfield $K\subseteq \mathfrak{C}$ is an $L_{\lambda,b}$-substructure of $\mathfrak{C}$ if and only if $b\subset K$ and $b$ is a $p$-basis of $K$ (in such a case the field extension $K\subseteq \mathfrak{C}$ is separable and even \'{e}tale). Each $L_{\lambda,b}$-substructure is definably closed and for $L_{\lambda,b}$-substructures, the model-theoretic algebraic closure coincides with the field theoretic separable closure (see \cite[Section 1.6]{ChHr}) and similarly with forking (see \cite[Section 1.8]{ChHr}).
Using the above, we immediately get the following description of regular extensions.
\begin{fact}\label{regext}
Let  $K_0\subseteq K_1$ be an extension of $L_{\lambda,b}$-substructures of $\mathfrak{C}$. Then, the extension $K_0\subseteq K_1$ is $\mathrm{SCF}_{p,e}$-regular if and only if $K_0\subseteq K_1$ is a regular extension of pure fields.
\end{fact}
We describe now PAC structures in the theory $\mathrm{SCF}_{p,e}$. This description appeared in \cite{Afshordel}, but the proof is only sketched there.
\begin{theorem}\label{efinite}
Let $K$ be a $L_{\lambda,b}$-substructure of $\mathfrak{C}$. Then, $K$ is $\mathrm{SCF}_{p,e}$-PAC if and only if $K$ is a PAC field.
\end{theorem}
\begin{proof}
In this proof, $K$ is a $\mathrm{SCF}_{p,e}$-substructure of $\mathfrak{C}$.
\\
$(\Rightarrow)$ Let us assume that $K$ is $\mathrm{SCF}_{p,e}$-PAC and let $K\subseteq N$ be a regular field extension. Then, there is a field extension $N\subset N'$ such that $b$ is a $p$-basis of $N'$ (we recall that now $b$ is a fixed $p$-basis of $\mathfrak{C}$). Therefore, we can assume that $K\subset N'$ is an extension of $L_{\lambda,b}$-substructures of $\mathfrak{C}$. By Fact \ref{regext}, $K\subseteq N'$ is also a $\mathrm{SCF}_{p,e}$-regular extension and $K$ is $L_{\lambda,b}$-existentially closed in $N'$. Therefore, $K$ is also is $L_{\lambda,b}$-existentially closed in $N$ and $K$ is existentially closed in $N$ in the field sense, hence $K$ is a PAC field.

$(\Leftarrow)$ Let us assume that $K$ is a PAC field. Let $b=(b_1,\ldots,b_e)$ be a $p$-basis of $K$, which is also a $p$-basis of $\mathfrak{C}$. We will consider the unary $\lambda$-functions
$$\lambda_{1,e},\ldots,\lambda_{p^e,e}:\mathfrak{C}\to \mathfrak{C}$$
(see \cite[Section 1.8]{Cha02}) with respect to this $p$-basis (they preserve $K$). We also take $a\in \mathfrak{C}^n$ such that $p(x):=\tp(a/K)$ is stationary and a quantifier-free $L_{\lambda,b}(K)$-formula $\phi(x)\in p(x)$. We need to show that $\phi$ has a realization in $K$. We inductively unravel all the terms appearing in the formula $\phi$. For example, if we have:
$$\phi(x)\colon \lambda_{i,n}\left(\lambda_{j,m}(x)+x^2\right)+x=0,$$
then we set:
$$\bar{a}:=\left(a,\lambda_{j,m}(a),\lambda_{i,n}\left(\lambda_{j,m}(a)+a^2\right)\right).$$
Then, for any
$$(b_1,b_2,b_3)\in \locus_K(\bar{a}),$$
we obtain that:
$$b_2=\lambda_{j,m}(b_1),\ \ \ \ b_3=\lambda_{i,n}\left(\lambda_{j,m}(b_1)+b_2^2\right).$$
As usual when we deal with fields, we can assume that there are only equalities in the formula $\phi$ (by replacing negations of equalities with equalities in ``higher dimensions''). Using the above procedure, we obtain a tuple $\bar{a}=(a_1,\ldots,a_t)$ such that $a=a_1$ and if we define:
$$V:=\locus_K(\bar{a}),$$
then for any $\bar{b}\in V(\mathfrak{C})$, we get that $\mathfrak{C}\models \phi(\bar{b})$. Since the type $p$ is stationary, $V$ is absolutely irreducible. Therefore, we obtain $\bar{b}\in V(K)$ such that $\mathfrak{C}\models \phi(\bar{b})$.
\end{proof}
There are three possible languages such that the theory $\mathrm{SCF}_{p,e}$ ($e$ finite) considered in each of these languages has quantifier elimination: $L_{\lambda,b}$, $L_{\lambda}$, and the Hasse-Schmidt language (see \cite{Zieg2}). One could wonder whether such a choice of the language affects the corresponding notion of a PAC-substructure. We address these issues in general below.
\begin{remark}
Assume that $T$ is a stable $L$-theory and $T'$ is an $L'$-theory being an extension by definitions of $T$.
Consider a model $M$ of $T$ and its counterpart $M'$ as an $L'$-structure (i.e. $M$ equipped with the natural $L'$-structure), and a subset $K$ of $M$.
\begin{enumerate}
\item If $K$ is a PAC substructure in the sense of $M$, then $K$ is an $L'$-substructure of $M'$ which is also a PAC substructure in the sense of $M'$.

\item If $K$ is a PAC substructure in the sense of $M'$, then $K$ is a PAC substructure in the sense of $M$.
\end{enumerate}
\end{remark}

For a finite group $G$, by Theorem \ref{thm:ec.if.PAC} we get existence of the theory $G-\mathrm{SCF}_{p,e}$, which is the model companion of the theory of actions of $G$ on characteristic $p$ fields of inseparability degree $e$. This theory was already analyzed in \cite{HKJSL} using different methods.

\subsubsection{Separably closed fields of infinite imperfection degree} .
We consider the theory SCF$_{p,\infty}$ in the language $L_{\lambda}$, where the $\lambda$-functions are multi-variable (see \cite[Section 1.4]{ChHr}), this definition is recalled below (we follow \cite[Section 1.8]{Cha02} here). For each $e>0$, we fix an enumeration $(m_{i,e})_{1\leqslant i \leqslant p^e}$ of the monomials $X_1^{i_1}\ldots X_e^{i_e}$ where $0  \leqslant i_1, \ldots , i_e  \leqslant p-1$, and define the functions $\lambda_{i,e}:K^e\times K \to K$ by considering the following three cases. Let $b_1,\ldots,b_e,c\in K$ and $1\leqslant i \leqslant p^n$.
\\
{\bf Case 1} \emph{$b_1,\ldots,b_e$ are $p$-dependent.}
\\
We set $\lambda_{i,e}(b_1,\ldots,b_e;c)=0$.
\\
{\bf Case 2} \emph{$b_1,\ldots,b_e,c$ are $p$-independent.}
\\
We set $\lambda_{i,e}(b_1,\ldots,b_e;c)=0$.
\\
{\bf Case 3} \emph{$b_1,\ldots,b_e$ are $p$-independent and $b_1,\ldots,b_e,c$ are $p$-dependent.}
\\
We use the following defining formula:
$$c=\sum_{j=1}^{p^e}\lambda_{j,e}(b_1,\ldots,b_e;c)^{p}m_{j,e}(b_1,\ldots,b_e).$$
Then, the theory SCF$_{p,\infty}$ has quantifier elimination in the language $L_{\lambda}$, but it does not have  elimination of imaginaries (see \cite[Section 1.8]{Cha02}). As for any theory of fields, SCF$_{p,\infty}$ eliminates finite imaginaries. Each $L_{\lambda}$-substructure is definably closed and for $L_{\lambda}$-substructures, the model-theoretic algebraic closure coincides with the field theoretic separable closure (see \cite[Section 1.6]{ChHr}) and similarly with forking (see \cite[Section 1.8]{ChHr}).

As in the previous part, we get the following description.
\begin{fact}\label{regext1}
Let $M\models \mathrm{SCF}_{p,\infty}$ and $K_0\subseteq K_1$ be an extension of $L_{\lambda}$-substructures of $M$. Then, $K_0\subseteq K_1$ is $\mathrm{SCF}_{p,\infty}$-regular if and only if $K_0\subseteq K_1$ is a regular extension of pure fields.
\end{fact}
We also need the following result of Tamagawa, which we phrase in geometric terms.
\begin{theorem}[Tamagawa, Proposition 11.4.1 in \cite{FrJa}]\label{tamagawa}
Let $V$ be an absolutely irreducible affine variety over a PAC field $K$ of characteristic $p>0$. Suppose that $f_1,\ldots,f_m\in K[V]$ are $p$-independent in $K(V)$ and $m$ is not greater than the imperfection degree of $K$. Then, there is $a\in V(K)$ such that $f_1(a),\ldots,f_m(a)$ are $p$-independent in $K$.
\end{theorem}
We need a slight enhancement of Tamagawa's Theorem, which we state and show below.
\begin{cor}\label{tamagawa2}
Let $V$ be an absolutely irreducible affine variety over a PAC field $K$ of characteristic $p>0$ and infinite imperfection degree. Suppose that we have a finite matrix $(f_{i,j})_{i,j}$ of elements of $K(V)$ whose rows are $p$-independent in $K(V)$. Then, there is $a\in V(K)$ such that each $f_{i,j}$ is defined at $a$ and the rows of the matrix $(f_{i,j}(a))_{i,j}$ are $p$-independent in $K$.
\end{cor}
\begin{proof}
We will often use the fact that the $p$-independence satisfies the Steinitz Exchange Principle and yields a pregeometry (see \cite[Remark C.1.1.3]{TentZie}), hence we have the corresponding dimension notions, which we denote by $\dim_p^{K}$ and $\dim_p^{K(V)}$.

We deal first with the case when $f_{i,j}\in K[V]$. For simplicity, let us assume that there are only three rows in our matrix and that each row has the same length $m$. Let us denote these rows by $(f_i)_i, (g_i)_i, (h_i)_i$. After permutation, there are $k,l\leqslant m$ such that:
\begin{enumerate}
  \item $f_1,\ldots,f_m,g_1,\ldots,g_k,h_1,\ldots,h_l$ are $p$-independent in $K(V)$;
  \item $f_1,g_1,\ldots,f_m,g_m\in \cl_p^{K(V)}(f_1,\ldots,f_m,g_1,\ldots,g_k)$;
  \item $f_1,g_1,h_1,\ldots,f_m,g_m,h_m\in \cl_p^{K(V)}(f_1,\ldots,f_m,g_1,\ldots,g_k,h_1,\ldots,h_l)$.
\end{enumerate}
In particular, we obtain that:
\begin{equation}
\dim_p^{K(V)}(f_1,g_1,h_1,\ldots,f_m,g_m,h_m)=m+k+l.\tag{i}
\end{equation}
By Theorem \ref{tamagawa} and Item $(1)$ above, there is $a\in V(K)$ such that
\begin{equation}
\text{$f_1(a),\ldots,f_m(a),g_1(a),\ldots,g_k(a),h_1(a),\ldots,h_l(a)$ are $p$-independent in $K$.}\tag{ii}
\end{equation}
Using (ii) and Item $(3)$ above, we obtain that:
\begin{equation}
\dim_p^{K}(f_1(a),g_1(a),h_1(a),\ldots,f_m(a),g_m(a),h_m(a))=m+k+l.\tag{iii}
\end{equation}
We will show that this choice of $a$ works. We focus on the most complicated case, that is we will prove that $h_1(a),\ldots,h_m(a)$ are $p$-independent in $K$ (the $p$-independence of $g_1(a),\ldots,g_m(a)$ is indeed easier by Item $(2)$ above). Since $h_1,\ldots,h_m$ are $p$-independent in $K(V)$, there are (after a permutation) $v\leqslant m,w\leqslant k$ such that $f_1,\ldots,f_v,g_1,\ldots,g_w,h_1,\ldots,h_m$ are $p$-independent in $K(V)$ and
\begin{equation}
f_1,g_1,h_1,\ldots,f_m,g_m,h_m\in \cl_p^{K(V)}(f_1,\ldots,f_v,g_1,\ldots,g_w,h_1,\ldots,h_m).\tag{iv}
\end{equation}
Using (i), (iii) and (iv), we obtain that:
$$v+w+m=\dim_p^{K}(f_1(a),g_1(a),h_1(a),\ldots,f_m(a),g_m(a),h_m(a)).$$
If $h_1(a),\ldots,h_m(a)$ were $p$-dependent in $K$, we would obtain by (iv) that:
$$\dim_p^{K}(f_1(a),g_1(a),h_1(a),\ldots,f_m(a),g_m(a),h_m(a))<v+w+m,$$
a contradiction.

We consider now that case when $f_{i,j}\in K(V)$. Let $U\subseteq V$ be an open $K$-subvariety which is $K$-isomorphic to an affine variety
such that is contained in the intersection of all $\dom(f_{i,j})$. Then $U$ is absolutely irreducible as well, and it is enough to apply the previously shown case of $f_{i,j}\in K[V]$.
\end{proof}

\begin{theorem}\label{einf}
$K$ be an $L_{\lambda}$-substructure of a monster model $\mathfrak{C}\models \mathrm{SCF}_{\infty,e}$. Then, $K$ is $\mathrm{SCF}_{p,\infty}$-PAC if and only if $K$ is PAC and $[K:K^p]=\infty$.
\end{theorem}

\begin{proof}
For the left-to-right implication, we take $K$ which is $\mathrm{SCF}_{p,\infty}$-PAC. By \cite[Section 1.7]{ChHr}, the properties of the generic 1-type in the theory $\mathrm{SCF}_{p,\infty}$ imply that $[K:K^p]=\infty$.
As in the previous part, we notice that extensions of $L_{\lambda}$-substructures of models of $\mathrm{SCF}_{p,\infty}$ are $\mathrm{SCF}_{p,\infty}$-regular if and only if they are a regular extension of pure fields. Now, the proof is identical to the proof of the corresponding implication in Theorem \ref{efinite}.

For the right-to-left implication, assume that  $K$ is PAC and $[K:K^p]=\infty$. Let us take $a\in \mathfrak{C}^m$ such that $p(x):=\tp(b/K)$ is stationary and a quantifier-free $L_{\lambda}$-formula $\phi(x)\in p(x)$ with parameters from $K$. We need to show that $\phi$ has a realization in $K$.
\\
{\bf Claim}
\\
There is $\bar{b}=(b,b')\in \mathfrak{C}^N$ such that $V:=\locus_K(\bar{b})$ is absolutely irreducible, and there is a finite matrix of rational functions $(f_{i,j}\in K(V))_{i,j}$ such that for each $i$, $f_{i,1},\ldots,f_{i,m_i}$ are $p$-independent in $K(V)$ and such that for all $\bar{c}=(c,c')\in V(\mathfrak{C})$, we have:
\\
IF for each $i$, $f_{i,1}(\bar{c}),\ldots,f_{i,m_i}(\bar{c})$ are $p$-independent in $\mathfrak{C}$, THEN $\mathfrak{C}\models \phi(c)$.
\begin{proof}[Proof of Claim]
By \cite[Lemma 2.9]{Cha02}, the formula $\phi(x)$ is equivalent in $\mathfrak{C}$ to an $L(K)$-formula of the form:
$$\exists y\ \alpha(x,y)\wedge \beta(x,y),$$
where $\alpha$ is quantifier-free in the language of fields and $\beta$ is a finite conjunction of universal formulas expressing that some subtuples of $xy$ are $p$-independent. Since $\phi(x)\in \tp(b/K)$, there is $b'\subset \mathfrak{C}$ such that:
$$\mathfrak{C}\models \alpha(b,b')\wedge \beta(b,b').$$
This is our choice of $b'$ as in the statement of this claim and the matrix of rational functions is given just by the coordinate functions expressing that $\beta$ is a ``conjunction of universal formulas expressing that some subtuples of $xy$ are $p$-independent'' (each row in this matrix corresponds to one formula from the finite conjunction giving the formula $\beta$).
\end{proof}
By Corollary \ref{tamagawa2}, there is $\bar{a}=(a,a')\in V(K)$ such that for each $i$, $f_{i,1}(\bar{a}),\ldots,f_{i,m_i}(\bar{a})$ are $p$-independent in $K$. Since the field extension $K\subseteq \mathfrak{C}$ is separable, each $f_{i,1}(\bar{a}),\ldots,f_{i,m_i}(\bar{a})$ is also $p$-independent in $\mathfrak{C}$. By Claim, we get that $\mathfrak{C}\models \phi(a)$, which we needed to show.
\end{proof}

\begin{theorem}
If $G$ is finite, then the model companion of the theory $\big((\scf_{p,\infty})_{\forall}\big)_G$,
denoted by $G-\scf_{p,\infty}$, exists.
\end{theorem}
\begin{proof}
We want to use Theorem \ref{thm:ec.if.PAC}, so we need that $\scf_{p,\infty}$
has QE, FS$+$, ST$+$ and that PAC is a first order property in $\scf_{p,\infty}$.
After Theorem \ref{einf}, the only thing which needs to be checked is ST$+$,
which might be a well-known fact, but as we could not find a proof of it, we noticed that one can adapt the proof of Lemma 4.16 from \cite{BHKK}.
\end{proof}

\subsubsection{Differentially closed fields}\label{dcfields}
We consider the theory $\mathrm{DCF}_{p}$  in the language $L_{\lambda_0,D}$, where $\lambda_0$ is the inverse of Frobenius on $p$-th powers and identically $0$ elsewhere.
Then, DCF$_{p}$ has quantifier elimination  \cite[Theorem 11]{Wood1},
but it does not have  elimination of imaginaries \cite[Remark 4.3]{MeWo}. It was shown in \cite{Shelahdif} that the theory $\mathrm{DCF}_{p}$ is stable.
\begin{remark}
The above ``$\lambda_0$-notation'' was introduced by the second author in \cite{K2} and perhaps it was not a very good choice, since:
\begin{itemize}
  \item it does not follow the original ``$r$-notation'' of Wood (see \cite[Section 2]{Wood2});
  \item $\lambda_{\emptyset}$ is defined as the identity function in \cite[Section 4]{messmer}.
\end{itemize}
However,  for the empty tuple $\bar{b}$ the most
natural interpretation of $\lambda_0$(=$\lambda_{0,0}$) is the one given above. Since this $\lambda_0$-notation was already used in several other papers, we stick with it in this paper as well.
\end{remark}
We think that the result below is a folklore, but we could not find a reference, so we give a proof instead.
\begin{fact}\label{dcffact}
Let $(\mathfrak{C},D)\models \mathrm{DCF}_{p}$ be a monster model and $K$ be an $L_{\lambda_0,D}$-substructure of $(\mathfrak{C},D)$. Then we have the following.
\begin{enumerate}
  \item The model-theoretic algebraic closure of $K$ coincides with its field theoretic separable closure.
  \item $K=\dcl(K)$.
\end{enumerate}
\end{fact}
\begin{proof}
Since $K$ is an $L_{\lambda_0,D}$-substructure of $(\mathfrak{C},D)\models \mathrm{DCF}_{p}$, the field extension $K\subseteq \mathfrak{C}$ is separable (see e.g. the beginning of the proof of \cite[Proposition 4.10]{BHKK}, where this separability appears in a much more general context). Therefore,
$K$ is a also a $L_{\lambda,D}$-substructure of $(\mathfrak{C},D)$, where $L_{\lambda,D}$ is the language with function symbols for all $\lambda$-functions.

For Item $(1)$, we note that the separable closure of $K$ is still a $L_{\lambda,D}$-substructure of $(\mathfrak{C},D)$. By (a more general) \cite[Lemma 4.14]{BHKK}, we get our description of the model-theoretic algebraic closure.

For Item $(2)$, if $a\in \dcl(K)$, then by Item $(1)$, we get that $a$ is separably algebraic over $K$. Since $D$ on $K$ extends uniquely to $K^{\sep}$, by quantifier elimination of $\mathrm{DCF}_p$ in $L_{\lambda_0,D}$ (or $L_{\lambda,D}$), we get that the type $\mathrm{tp}^{\mathrm{DCF}_p}(a/K)$ is isolated by $f_a$: the minimal polynomial of $a$ over $K$. Since $a\in \dcl(K)$, we get that $\deg(f_a)=1$, so $a\in K$, which we needed to show.
\end{proof}
Again, we need the following description of regular extensions with respect to the theory we consider. It follows immediately from Fact \ref{dcffact}.
\begin{fact}\label{regextdif}
Let $(\mathfrak{C},D)\models \mathrm{DCF}_{p}$ be a monster model and $(K_0,D)\subseteq (K_1,D)$ be an extension of $L_{\lambda_0,D}$-substructures of $\mathfrak{C}$. Then, $K_0\subseteq K_1$ is $\mathrm{DCF}_{p}$-regular if and only if $K_0\subseteq K_1$ is a regular extension of pure fields.
\end{fact}

We specify now an $L_{\lambda_0,D}$-theory of some differential fields in positive characteristic. We need the following working definition first.
\begin{definition}\label{adm}
Let $K$ be a field of characteristic $p>0$. A tuple $(V;f_1,\ldots,f_n)$ is \emph{admissible}, if $V$ is a $K$-irreducible affine $K$-variety and $f_1,\ldots,f_n\in K(V)\setminus K(V)^p$.
\end{definition}
We note the following obvious property.
\begin{remark}\label{remadm}
Let $K\subseteq M$ be a field extension. For any $f\in K(V)$ and any $a\in V(M)$ generic of $V$ over $K$, we have that $f\in K(V)^p$ if and only if $f(a)\in K(a)^p$.
\end{remark}
\begin{lemma}\label{admfo}
Assume that $K$ is PAC of infinite imperfection degree (actually, non-perfect would be enough). Then, the above notion of an admissible tuple is first-order in parameters of this tuple.
\end{lemma}
\begin{proof}
Let us take $f\in K(V)$. By Corollary \ref{tamagawa2}, we obtain that $f\in K(V)^p$ if and only if $f(V(K))\subseteq K^p$. Since the second condition is clearly first-order, the result follows.
\end{proof}
The next question is not related to our results and we find it a bit amusing. The answer may be simple, but we were unable to find it.
\begin{question}
Is the property ``$f\in K(V)^p$'' first-order in parameters of $f$ and $V$ for an algebraically closed $K$ (that is: modulo the theory ACF$_p$)?
\end{question}
To state our axioms for PAC-$\mathrm{DCF}_p$ differential fields, we need to recall some notions. We decided to work here with the case of differential fields for the clarity of presentation, however, as we will see in Section \ref{secfop}, these results hold in a much greater generality. Still, our references here come from this more general context, since we do not know any source where they are stated exactly for the differential case.

Let $(K,D)$ be a differential field (for a while the characteristic of $K$ does not matter) and $V$ be a $K$-variety. Then, $\tau^D(V)$ denotes the \emph{prolongation} of $V$ with respect to $D$, which in this case can be described as a torsor of the tangent bundle of $V$ (see \cite[Definition 1.4]{PP} and \cite[Definition 4.1]{MS1}). We have a natural map (see e.g. \cite[Remark 2.13]{BHKK}):
$$D_V:V(K)\ra \tau^D(V)(K).$$
Let $K\subseteq \Omega$ be a field extension and $a,a'\subset \Omega$ be such that:
$$V=\locus_K(a),\ \ \ W=\locus_K(a,a').$$
For reader's convenience, we recall now two results from \cite{BHKK} which we will use.
\begin{lemma}[Lemma 3.5 in \cite{BHKK}]\label{easylemma2}
The following are equivalent.
\begin{enumerate}
\item There is a derivation $D':K[a]\subseteq K[a,a']$ extending $D$ such that $D'(a)=a'$.

\item $W\subseteq \tau^{D}(V)$.
\end{enumerate}
\end{lemma}
Assume that $V,W$ are $K$-varieties as in the statement of Lemma \ref{easylemma2}, that is: $W\subseteq \tau^{D}(V)$. Let $\iota:W\to \tau^{D}(V)$ denote the inclusion morphism and
$$\alpha:=\pi^V_{D}\circ \iota:W\to V.$$
Consider the following (not necessarily commutative!) diagram:
\begin{equation*}
 \xymatrix{  &  &  \tau^{D}(W)  \ar[lld]_{\pi^W_{D}} \ar[rrd]^{\tau^{D}(\alpha)} &  &   \\
W\ar[rrrr]^{\iota}  &  & &  &  \tau^{D}(V) .}
\end{equation*}
Using this diagram, we define the following $K$-subvariety of $\tau^{D}(W)$:
\begin{IEEEeqnarray*}{rCl}
E &:=& \mathrm{Equalizer}\left(\tau^{D}(\alpha),\iota \circ \pi^W_{D}\right)\\
  &=& \left\{a\in \tau^{D}(W)\ |\ \tau^{D}(\alpha)(a)=\iota \circ \pi^W_{D}(a)\right\}.
\end{IEEEeqnarray*}
The following result is crucial.
\begin{theorem}[Proposition 3.6 in \cite{BHKK} specialized to the case of derivations]\label{kerprol}
The following are equivalent.
\begin{enumerate}
\item The morphism $\pi_E:E\to W$ is dominant.

\item There is a derivation on $K(a,D(a))$ extending $D:K[a]\to K[a,D(a)]$.
\end{enumerate}
\end{theorem}
We state our axioms below.
\smallskip
\\
\textbf{Axioms for $D-PAC$}
\\
Let $(K,D)$ be a differential field of characteristic $p>0$ and for each pair of affine $K$-varieties $(V,W)$ and each tuple $f_1,\ldots,f_n\in K(V)$ such that
\begin{itemize}
\item $W$ is absolutely irreducible,

\item $W\subseteq \tau^{D}(V)$,

\item the projection $\pi:W\to V$ is dominant,

\item $E$ projects dominantly on $W$,

\item the tuple $(W;f_1\circ \pi,\ldots,f_n\circ \pi)$ is admissible;
\end{itemize}
there is $x\in V(K)$ such that $f_1(x),\ldots,f_k(x)$ are not $p$-th powers in $K$ and $D_V(x)\in W(K)$.
\smallskip
\\
By standard arguments (see e.g. \cite[Remark 2.7(1)]{HK3}) and Lemma \ref{admfo}, the above axiom scheme is first-order. An essential argument using Theorem \ref{kerprol} in the proof of the result below follows the ideas of the proof of  \cite[Proposition 5.6]{PilPol}.   
\begin{theorem}\label{dpac}
Let $(K,D)$ be a differential field of characteristic $p>0$ considered as an $L_{\lambda_0,D}$-structure. Then, $(K,D)$ is $\mathrm{DCF}_p$-PAC if and only if $(K,D)$ is $D$-PAC (as defined above).
\end{theorem}
\begin{proof}
For the left-to-right implication, we assume that $(K,D)$ is $\mathrm{DCF}_{p}$-PAC and $((V;f_1,\ldots,f_n),W)$ is as in the assumptions of the axioms of $D$-PAC. By Lemma \ref{easylemma2}, there is a derivation $D'\colon K(V)\to K(W)$ of the inclusion $K(V)\subseteq K(W)$ (given by the dominant morphism $\pi:W\to V$). By Theorem \ref{kerprol}, $D'$ extends to a derivation $D''\colon K(W)\to K(W)$. Since $W$ is absolutely irreducible, the extension $K\subseteq K(W)$ is regular and $(K,D)\subseteq (K(W),D'')$ is also a differential field extension. By Fact \ref{regextdif}, $(K,D)\subseteq (K(W),D'')$ is a $\mathrm{DCF}_p$-regular extension. Since $(K,D)$ is $\mathrm{DCF}_p$-PAC, we get that $(K,D)$ is existentially closed in $(K(W),D'')$ (in the language $L_{\lambda_0,D}$).

Let us choose:
$$a=\id_{K[V]}\in V(K(V))\subseteq V(K(W)).$$
Then, as usual, we have $D''_V(a)\in W(K(W))$. Since $a$ is a generic point of $V$ over $K$, by Remark \ref{remadm}, we get that $f_1(a),\ldots,f_k(a)$ are not $p$-th powers in $K(W)$. Since $(K,D)$ is existentially closed in $(K(W),D'')$ (in the language $L_{\lambda_0,D}$), there is $\alpha\in V(K)$ such that $f_1(\alpha),\ldots,f_k(\alpha)$ are not $p$-th powers in $K$
and $D_V(\alpha)\in W(K)$.
\\
\\
For the right-to-left implication, we assume that  $(K,D)$ is $D$-PAC. Let us take $a\in M^n$ such that $p(x):=\tp(a/K)$ is stationary and a quantifier-free $L_{\lambda_0,D}(K)$-formula $\phi(x)\in p(x)$. We need to show that $\phi$ has a realization in $K$. As usual, we can assume that $\phi$ does not contain negations of equalities. We will ``correct'' now the formula $\phi(x)$ (at the cost of adding extra variables, some fixed terms, and a new tuple of elements of $M$ including $a$) into a new formula $\varphi(\bar{x})\in L_D$ over $K$ such that $\bar{x}=(x_1,\ldots,x_l)$ and $x_1=x$.

We illustrate this ``correction'' using an example first. Assume that the formula $\phi(x)$ has the following form:
$$D\left[\lambda_0\left(D\left(\lambda_0(x)\right)+D(x)\right)\right]+x=0.$$
We consider two cases, where each of them has two subcases. For the new variables, we will use $y,z$ rather than $x_2,x_3$.
\\
\textbf{Case I}: $\lambda_0(a)=0$.
\\
\textbf{Subcase I.1}: $\lambda_0(D(a))=0$ (so: $a=0$).
\\
The ``correction'' is $\varphi(x)\colon x=0$  and there are no extra variables and no fixed terms.
\\
\textbf{Subcase I.2}: $\lambda_0(D(a))\neq 0$.
\\
The ``correction'' is
$$\varphi(x,y)\colon y^p=D(x) \wedge D(y)+x=0,$$
the fixed term is $t_1(\bar{x})=x$, and $\bar{a}=(a,D(a)^{1/p})$.
\\
\textbf{Case II}: $\lambda_0(a)\neq 0$ (so: $D(a)=0$).
\\
\textbf{Subcase II.1}: $\lambda_0\left(D\left(\lambda_0(a)\right)\right)=0$.
\\
It cannot happen, since then $a=0$ and $\lambda_0(a)\neq 0$.
\\
\textbf{Subcase II.2}: $\lambda_0\left(D\left(\lambda_0(a)\right)\right)\neq 0$.
\\
The ``correction'' is
$$\varphi(x,y,z)\colon y^p=x \wedge z^p=D(y)+D(x) \wedge D(z)+x=0,$$
there are no fixed terms, and
$$\bar{a}=\left(a,a^{1/p},\left(D\left(a^{1/p}\right)\right)^{1/p}\right).$$
The general procedure can be explained using induction on the complexity of $L_{\lambda_0,D}$-terms over $K$. We obtain a quantifier-free $L_D$ formula $\varphi(\bar{x})$ over $K$, $L_D$-terms $t_1(\bar{x}),\ldots,t_k(\bar{x})$ over $K$ such that:
$$(*)\ \ \ \ \ \  \text{if $(K,D)\models \varphi(\bar{\alpha})$ and $t_1(\bar{\alpha})\notin K^p,\ldots,t_k(\bar{\alpha})\notin K^p$, then $K\models \phi(\alpha)$},$$
and $\bar{a}$ such that $(M,D)\models \varphi(\bar{a})$ and  
$t_1(\bar{a})\notin M^p,\ldots,t_k(\bar{a})\notin M^p$.

Let us take now  a quantifier-free $L$ formula $\psi(\tilde{x})$ over $K$ such that:
$$\varphi(\bar{x})\colon \ \ \ \psi\left(\bar{x},D(\bar{x}),\ldots,D^m(\bar{x})\right)$$
for some $m\in \Nn$. Let us define:
$$\tilde{a}:=\left(\bar{a},D(\bar{a}),\ldots,D^m(\bar{a})\right),\ \ \ \ V:=\locus_K\left(\tilde{a}\right),\ \ \ \
W:=\locus_K\left(\tilde{a},D(\tilde{a})\right).$$
Let $\pi:W\to V$ denote the dominant projection on the ``$\tilde{a}$-coordinates''. There are rational function symbols $f_1(\tilde{x}),\ldots,f_k(\tilde{x})$ over $K$ such that for each $i$, we have:
$$t_i(\bar{x})=f_i\left(\bar{x},D(\bar{x}),\ldots,D^m(\bar{x})\right).$$
Therefore, we obtain that for each $i$:
$$f_i\left(\tilde{a}\right)=t_i\left(\bar{a}\right)\notin M^p\supseteq K\left(\tilde{a},D(\tilde{a})\right)^p=K(W)^p,$$
so $(W;f_1\circ \pi,\ldots,f_n\circ \pi)$ is an admissible tuple.

Let us take $\tilde{\alpha}\in V(K)$ such that $f_1(\tilde{\alpha}),\ldots,f_k(\tilde{\alpha})$ are not $p$-th powers in $K$ and $D_V(x)\in W(K)$. By the construction we get that:
$$\tilde{\alpha}=\left(\bar{\alpha},D(\bar{\alpha}),\ldots,D^m(\bar{\alpha})\right).$$
Therefore, we obtain $(K,D)\models \varphi(\bar{\alpha})$ and $t_1(\bar{\alpha}),\ldots,t_k(\bar{\alpha})$ are not $p$-th powers in $K$. By $(*)$ above, we obtain that $K\models \phi(\alpha)$.
\end{proof}
\begin{remark}
It was shown in \cite{Gogolok} that the ``equalizer condition'' on the dominant map $E\to W$ can be replaced with the easier condition of separability of the map $W\to V$ and then we still get geometric axioms of $\mathrm{DCF}_p$ (it also applies to the case of derivations of the Frobenius map). However, we do not know whether such a replacement would also work for the PAC-axioms, since we will not have Theorem \ref{kerprol} after such a replacement.
\end{remark}

\begin{theorem}
If $G$ is finite, then the model companion of $\big( (\dcf_p)_{\forall}\big)_G$, denoted by $G-\dcf_p$, exists.
\end{theorem}

\begin{proof}
Once again, we would like to use Theorem \ref{thm:ec.if.PAC}.
We have that $\mathrm{DCF}_{p}$ satisfies FS$+$ and it will be shown in a greater generality (see again Section \ref{secfop}) that types over algebraically closed sets are stationary (so ST$+$ follows). Theorem \ref{dpac} assures us that the PAC property is first order in $\dcf_p$.
\end{proof}
We would like to include here the following general result which will be immediately useful.
\begin{remark}\label{eq}
Let $T$ be a $L$-theory with quantifier elimination and $G$ be an arbitrary group.
\begin{enumerate}
  \item If the theory $G-T$ exists, then the theory $G-(T^{\eq})^m$ exists as well, where the superscript ``$m$'' denotes the Morleyization.

  \item If $G$ is finite, $T$ is stable, and the theory $G-T$ exists, then the theory $G-(T^{\eq})^m$ is simple.
\end{enumerate}
\end{remark}
\begin{proof}
The proof of Item $(1)$ is straightforward and we leave it to the reader. Item $(2)$ follows from \cite[Cor 4.28]{Hoff3}.
\end{proof}
Using Remark \ref{eq}, we obtain the following.
\begin{cor}
Let $G$ be a finite group and $p$ be a prime number.
\begin{enumerate}
\item The theory $G-\scf_{p,\infty}$ is strictly simple, that is simple, not stable, and not supersimple.

\item The theory $G-\dcf_p$ is strictly simple as well.
\end{enumerate}
\end{cor}

\subsection{Fields with operators}\label{secfop}
In this subsection, we briefly explain how to generalize Theorem \ref{dpac} beyond the case of differential fields. We recall below some of the set-up from \cite{BHKK}.

Let $\ka$ be a field and $B$ be a finite local $\ka$-algebra of dimension $e$. Assume that we have a $\ka$-algebra map $\pi_B:B\to \ka$. Let $\{b_0,\ldots,b_{e-1}\}$ be a fixed $\ka$-basis of $B$ such that $b_0=1$ and $\pi_B(b_i)=0$ for $i>0$. For convenience, we also set $d:=e-1$.
\begin{definition}\label{bopdef}
Assume that $R$ and $T$ are $\ka$-algebras  and let $\partial=(\partial_0,\ldots,\partial_d)$ where $\partial_0,\ldots,\partial_d:R\to T$ are $\ka$-algebra homomorphisms.
\begin{enumerate}
\item If $R=T$ and $\partial_0=\id$, then we say that $\partial$ is a \emph{$B$-operator on $R$} if the corresponding map
$$R\ni r\mapsto \partial_0(r)\otimes b_0+\ldots+\partial_d(r)\otimes b_d\in R\otimes_{\ka}B$$
is a $\ka$-algebra homomorphism. We will also denote the map above by the same symbol $\partial$.

\item More generally, if the corresponding map
$$R\ni r\mapsto \partial_0(r)\otimes b_0+\ldots+\partial_d(r)\otimes b_d\in T\otimes_{\ka}B$$
is a $\ka$-algebra homomorphism, then we say that $\partial$ is a \emph{$B$-operator from $R$ to $T$}.
Note that if $\partial$ is a $B$-operator from $R$ to $T$, then $\partial_0:R\to T$ is a $\ka$-algebra homomorphism.
\end{enumerate}
\end{definition}
Assume that $(K,\partial)$ is a field with a $B$-operator and $V$ is an affine $K$-variety. The notion of a prolongation $\tau^{\partial}(V)$ was defined in this generality (see \cite{MS2}) Under the additional assumption of $\mathrm{Fr}_B(\ker(\pi_B))=0$ (see \cite[Remark 3.3]{BHKK}), we get  the versions of Lemma \ref{easylemma2} and Theorem \ref{kerprol} (actually, the references from Section \ref{dcfields} are coming exactly from the $B$-operator context). We get the corresponding Axioms for $\partial$-PAC (with the identical formulation) and a generalization of Theorem \ref{dpac}. The proof of this generalization is conceptually the same, but would be more cumbersome to write comparing to the proof of Theorem \ref{dpac}. Therefore we decided to include the general case of $B$-operators as this comment only.
\begin{remark}
The argument above works in the even more general case of $\mathcal{B}$-operators (replacing $B$-operators from Definition \ref{bopdef}) as considered in \cite{GogK}. The main example of a $\mathcal{B}$-operator which is not a $B$-operator is a derivation of the Frobenius map.
\end{remark}
\begin{remark}
Using Theorem \ref{thm:ec.if.PAC}, we obtain that for a finite group $G$, the theory of $G$-actions on fields with $B$-operators (and also on fields with $\mathcal{B}$-operators) has a model companion. It may be shown more directly and without going through $\lambda$-functions. The axiomatization is as follows.
\smallskip
\\
\textbf{Axioms for $G$-$B$-DCF}
\\
The structure $(K,\partial,\sigma)$ is a $G$-$B$-field such that for each pair $(V,W)$ of $K^G$-varieties, IF
\begin{itemize}
\item the action of $G$ on $K$ is faithful,

\item $V$ and $W$ are $K$-irreducible,

\item $W\subseteq \tau^{\partial}(V)$,

\item $W$ projects generically on $V$,

\item $E$ projects generically on $W$;
\end{itemize}
THEN there is $x\in V(K^G)$ such that $\partial_V(x)\in W(K^G)$.
\end{remark}
\subsection{Other examples and questions}\label{secoeq}
We discuss now some other examples and ask some questions. 
In \cite{PilPol}, an example of a stable theory $T$ is given such that the class of $T$-PAC structures is not elementary. However, the theory $T$ in this example does not have quantifier elimination, so from our perspective it is not a good theory to test whether the PAC property is first-order.

More precisely, in \cite[Example 5.1]{PilPol} the theory of an equivalence relation with exactly one finite class of $n$ elements for each $n>0$ appears. This theory is (implicitly) considered in the natural language with one unary relation symbol. Then, this theory is not even model complete, since finite classes in a model may become infinite in its extension. Similarly, one can see that this theory is not inductive. Therefore, to have any hopes for quantifier elimination, one needs to add to the language the unary predicates $(R_n)_{n>0}$ naming all finite equivalence classes (see also \cite[Example 5.2]{PilPol}, where each \emph{element} of each finite equivalence class is named). Using Robinson's Test, it is not difficult to check that with such a choice of the language this theory becomes model complete and also substructure complete, so it has quantifier elimination. Then, taking algebraic closure is the same as adding ``missing points'' in all named finite classes. Therefore, if $M\subseteq M'$ and $M'$ is a subset of a model, then $\acl(M)\cap M'=M$ if and only if for all $n$, we have $R_n(M)=R_n(M')$. Thus, $M$ is PAC if and only if $M=\dcl(M)$ and $M$ is infinite, so PAC is a first-order property in this case (note that definably closed substructures are such ones that each finite class does not have ``co-size one'').

As a conclusion, we do not know any stable theory $T$ with quantifier elimination such that the class of $T$-PAC structures is not elementary. We formulate the relevant question below.
\begin{question}\label{q1}
Assume that $T$ is stable and has quantifier elimination. Is the class of $T$-PAC structures elementary?
\end{question}
Positive answer to the question above implies (using Theorem \ref{thm:ec.if.PAC}) that for a finite group
$G$ and for $T$ as above eliminating strong types and coding finite sets, the theory of actions of $G$ on models of $T_{\forall}$ has a model companion. 

\begin{remark}
This is related to a general conjecture of the first author, see \cite[Conjecture 5.2]{Hoff3}:
\\
``Assume that $T_0$ is theory with a model companion and $G$ is a finite group. Does the theory of $G$-actions on models of $T_0$ have a model companion?'' 
\\
This conjecture was meanwhile refuted in \cite[Remark 3.9(2)]{BK3}, where $T_0$ is the theory of difference fields. In this example, the model companion of $T_0$ (the theory ACFA) is neither stable nor it has quantifier elimination.
\end{remark}

To pursue the answer for the above question one can start with somehow stronger assumption:

\begin{question}\label{q1b}
Assume that $T$ has nfcp (no finite cover property) and has quantifier elimination. Is the class of $T$-PAC structures elementary?
\end{question}
\noindent
The above assumption on not having the finite cover property is related to the PAC property a little bit in Remark 3.6 in \cite{PilPol}, but the main point here is that it was shown in general that a stronger variant of the notion of nfcp (i.e. $T$ does not admit obstructions)
implies the model companion of the theory of models of $T$ with a group action of $\mathbb{Z}$ exists (\cite{balshe}).

Let $T$ be a stable theory with quantifier elimination. If we replace a finite group $G$ with the cyclic infinite group $\Zz$, then the model theory of actions of $\Zz$ on models of $T$ (we do not have distinguish between $T_{\forall}$ and $T$ in this case) has been thoroughly studied (see e.g. \cite{acfa1} and \cite{ChPi}). An analogue of our Question \ref{Q324} was asked in before Lemma 4.2 in \cite{PilPol}, that is it is asked there whether the existence of the theory $TA$ (which is called $\Zz-T$ in this paper) implies that $T$-PAC is first order. The main result of this paper, Theorem \ref{thm:ec.if.PAC}, gives the opposite implication in the case of finite groups. Such an implication is not true in the case of the actions of $\Zz$ (see \cite[Example 5.2]{PilPol}).
\begin{remark}\label{last}
To keep this paper reasonably sized, we have not checked all the known stable theories. However, there is one theory which we would not mind to analyze but the methods of this paper do not suffice to do that. This is the theory DCF$_{p,m}$ for $m>1$, that is the theory of differentially closed fields of characteristic $p>0$ with $m$ commuting derivations, see \cite{Pierce3} where the theory of differentially closed fields with $m$ commuting derivations is considered in arbitrary characteristic (it is called $m$-DCF in \cite{Pierce3}). Similarly as in the case of DCF$_{0,m}$ for $m>1$, we can not repeat our argument from the proof of Theorem \ref{dpac}, since we do not have a version of Theorem \ref{kerprol} in the case of several commuting derivations. It looks natural here to apply the approach of \cite{SanchezTressl} described briefly in Remark \ref{noethlast}(3) that is:
\begin{itemize}
  \item develop the appropriate notion of largeness for differential fields of positive characte\-ristic;
  \item show that the above notion is first-order;
  \item show that DCF$_{p,m}$-PAC is the same as the largeness above together with the PAC in the sense of fields.
\end{itemize}
We plan to pick it up in a further research.
\end{remark}

\bibliographystyle{plain}
\bibliography{1nacfa55}

\end{document}